\documentclass[12pt,psfig]{article}
 \input{psfig}
 \usepackage{amssymb}
 \newcommand{\R}{{\mathbb R}}
 \newcommand{\C}{{\mathbb C}}
 \newcommand{\Z}{{\mathbb Z}}
 \newtheorem{theorem}{Theorem}[section] 
 \newtheorem{lemma}[theorem]{Lemma}
 \newtheorem{proposition}[theorem]{Proposition}
 \newtheorem{corollary}[theorem]{Corollary} 
 \newtheorem{remark}[theorem]{Remark}

 \newtheorem{definition}[theorem]{Definition}

 \newtheorem{fact}[theorem]{Fact}
 \def\Box
  {\thinspace\vbox{\hrule height .5pt \hbox{\vrule  
   width .5pt \vbox to 7pt{\hbox to 3.5pt{}} \vrule width .5pt} 
   \hrule height 0pt depth .5pt}}

 \newenvironment{proof}{{\it Proof:\/}}{$\Box$\vskip 0.08in}
 \newcommand{\mod}{{\mbox{ mod }}}
\newcommand{\card}[1]{\overline{\overline{#1}}}
 \newcommand{\knot}[1]{\parbox{1.2cm}{\psfig{figure=cfig#1.eps,height=1.2cm}}}
 \newcommand{\rysunek}[2]{\psfig{figure=charfig#1.eps,height=#2}}

\begin{document}
\centerline{\LARGE\bf On Skein Algebras And $Sl_2(\C)$-Character}
\centerline{\LARGE\bf Varieties}
\vspace*{0.2in}
\centerline{J\'ozef H.~Przytycki and Adam S.~Sikora\footnote{
The second author was partially supported by NSF grant DMS93-22675.}}
\vspace*{0.2in}

\section{Introduction}
Let $M$ be an oriented $3$-manifold.
For any commutative ring $R$ with a specified invertible element $A$ one 
can assign an $R$-module ${\cal S}_{2,\infty}(M;R,A)$ called the Kauffman 
bracket skein module of $M.$ This invariant of $3$-manifolds was 
introduced by the first author in [P-1]. 

This paper gives insight into broad and intriguing connections
between two apparently unrelated theories: the theory of skein modules
of $3$-manifolds and the theory of representations of groups into special 
linear groups of $2\times 2$ matrices. This connection was first observed
by D. Bullock ([B-1],[B-2],[B-4]).

We believe that our research, which was originated in [P-S-1] and 
continued in this paper, will ultimately result in a theory which will 
reveal some of the mysterious mutual correlations between a skein 
approach and an approach via methods of representation theory to 
3-dimensional topology. 

Such a theory is needed, for example, in order to advance the study of 
quantum invariants. These invariants can be defined both in terms of skein 
theory and representation theory (of quantum groups). A lack of good 
understanding of relations between these two theories gives rise to some 
difficulties in studying quantum invariants.

Our work can be also considered in a context of the theory of Culler and 
Shalen ([C-S]) which relates properties of the 
$Sl_2(\C)$-character variety of $\pi_1(M),$ for a given 3-manifold $M,$ with 
properties of incompressible surfaces in $M.$ In this paper we give a 
topological interpretation of $Sl_2(\C)$-character variety of $\pi_1(M).$ 
Using this interpretation one can restate some of the deep results of 
Culler and Shalen in a purely topological manner.

We were inspired to write this paper after we discovered the results
of Brumfiel and Hilden, [B-H], which interplay very nicely with our own
work [P-S-1]. This paper is based on these two papers. In particular, we use 
the following results from [P-S-1]:
\begin{enumerate}
\item If $A=-1$ then ${\cal S}_{2,\infty}(M;R,A)$ is a commutative 
$R$-algebra, called a skein algebra of $M.$
\item For any group $G$ and a ring $R$ one can define a skein algebra of
$G,\ {\cal S}(G;R),$ in such a way that if $G=\pi_1(M)$ then the two notions 
of skein algebras coincide, i.e. 
$${\cal S}_{2,\infty}(M;R,-1)\simeq {\cal S}(G;R).$$
\end{enumerate}

Skein algebras will be the main subject of our study.
Here is the plan of our paper:\\

In Section 2 we give the definition of the Kauffman bracket skein module,
${\cal S}(M; R,A),$ and the skein algebra, ${\cal S}(G; R).$
We recall also all necessary results from [P-S-1] needed in this paper.
 
In Section 3 we prove that if some general conditions are satisfied
(for example, if ${1\over 2}\in R$) then ${\cal S}(G;R)$ is 
isomorphic to an algebra $TH_R(G)$ introduced by Brumfiel and 
Hilden in connection with their study of $Sl_2(R)$-representations of groups,
[B-H]. Using this fact we prove a few new results about the algebras $TH_R(G).$
For example, $TH_R(\pi_1(F))$ is a free module for any ring $R$ and for any
surface $F.$

In Section 7 we consider characters of $Sl_2$-representations of groups.
We show, using the work of Brumfiel and Hilden, 
that if $K$ is an algebraically closed field of characteristic $0$ then the 
algebra ${\cal S}(G;K)$ is isomorphic to the (unreduced) coordinate
ring of the $Sl_2(K)$-character variety of $G.$ Moreover, 
using various algebraic and topological methods, we prove that for many 
important classes of groups this algebra is actually reduced, i.e. 
the algebra ${\cal S}(G;K)$ does not have any nilpotent elements.

In Section 4 we consider another structure of a generally non-commutative
algebra assigned to ${\cal S}_{2,\infty}(F\times I;R,A),$ when $F$ is a 
surface. We present some of its interesting properties. For example,
${\cal S}_{2,\infty}(F\times I;R,A)$ is a central algebra over a ring 
of polynomials induced by the boundary of $F.$ Moreover,
${\cal S}_{2,\infty}(F\times I;R,A)$ has no zero divisors. In particular, this
implies that the $Sl_2(C)$-character varieties of the fundamental groups
of surfaces are irreducible.

In Sections 6 and 8 we find, using results of previous sections and 
some arguments from algebraic geometry, the minimal numbers of generators of 
various skein algebras. In particular, we prove that the minimal number of 
generators of ${\cal S}(F_n,R)$ is $n+{n \choose 2}+{n \choose 3}$ or $2^n-1$ 
depending on whether $2$ is invertible in $R$ or not.
We also give the minimal numbers of generators of skein algebras
${\cal S}_{2,\infty}(F\times I;R,A),$ for $R=\Z[A^{\pm 1}],$ and
estimate the minimal numbers of generators for other rings of coefficients, 
$R.$

In Section 5 we introduce for any surface $F$ a new algebra 
${\cal S}_{2,\infty}^{rel}(F\times I;R,A)$ which is built on relative links in
the manifold $M=F\times [0,1].$ We prove that this algebra is isomorphic to 
an algebra $H_R(G),\ G=\pi_1(M),$ investigated in the book by Brumfiel
and Hilden, [B-H]. This result, combined with our result from Section
3, gives a nice topological interpretation of algebraic objects,
$H_R(G)$ and $TH_R(G),$ considered in [B-H].

\section{Skein Modules and Skein Algebras}

In this paper we make the following assumptions:
\begin{enumerate}
\item All rings are commutative and have identities. All homomorphisms 
between rings preserve identities.
\item For a given ring $R$ and a set ${\cal S}$ we denote
the free $R$-module with a basis composed of elements of ${\cal S}$
by $R{\cal S}.$
\item All topological spaces (manifolds, links) are considered as objects in
the category of piecewise-linear topological spaces. In particular,
all links are tame. All continuous functions preserve PL-structures.
\item All manifolds are oriented and they may have boundaries.
\item Let $N$ be a submodule of $M.$ Then for any $m\in M$ we denote
$m+N\in M/N$ by $[m].$ Similarly, if $I$ is an ideal in a ring $R$ then
$x +I\in R/I$ is denoted by $[x],$ for any $x\in R.$
\end{enumerate}

We start with a definition of the Kauffman bracket skein module.
Skein modules were independently introduced by V. Turaev in [T-1] and,
in a more general setting, by the first author in [P-1]. Kauffman bracket
skein modules of manifolds were for the first time defined in [P-1].
They are composed of formal linear combinations of framed unoriented
links considered up to some local relations called
{\it skein relations}. By a framed link in a $3$-manifold $M$ we mean
an embedding of a finite family of annuli into the interior of $M.$
In our definition of the Kauffman bracket skein module we follow [H-P-1].

\begin{definition}\label{2.1}\ \\ 
Let $M$ be any oriented 3-manifold and let ${\cal L}_{fr}(M)$ denote the set 
of all ambient isotopy classes of framed unoriented links in the
interior of $M,$ including 
the empty link, $\emptyset$. Let $A$ be a specified invertible element in
a ring $R.$ Furthermore, let $S_{2,\infty}$ be a submodule of 
$R{\cal L}_{fr}(M)$ generated by two kinds of elements:
\begin{enumerate}
\item skein expressions $L_+-AL_0 - A^{-1}L_{\infty},$ 
where $L_+, L_0, L_{\infty}$ are any three framed links in $M$ which are the 
same outside a small 3-ball $B\subset M,$ but inside $B$ they are as
in Fig. 1,
i.e. there is an orientation preserving homeomorphism between $B$ and a ball
$B'\subset R^3$ which carries $B\cap L_+, B\cap L_-, B\cap L_0$ to
fragments of links presented in Fig. 1. The orientation of $B'\subset
\R^3$ is given by an ordered basis of $\R^3,$ $v_1,v_2,v_3,$ where 
$v_1, v_2$ are pictured in Fig. 1 and $v_3$ is orthogonal to $v_1$ and
$v_2$ and it is directed towards the reader. In this situation we say
that $L_+, L_0, L_{\infty}$ are skein related.
\item $L\cup \bigcirc +(A^2+A^{-2})\cdot L$, where $L$ is any link in $M$ 
and $L\cup \bigcirc$ denotes a disjoint union of $L$ with a trivial 
component, $\bigcirc$. We assume that $\bigcirc$ has a trivial 
framing and $\bigcirc$ is unlinked with $L$ i.e. $\bigcirc$ lies in a
plane in a ball disjoint from $L.$
\end{enumerate}

We define the Kauffman bracket skein module of $M$ as a
quotient \[{\cal S}_{2,\infty}(M;R,A)= R{\cal L}_{fr}(M)/S_{2,\infty}.\]
\end{definition}

\begin{figure}
\centerline{\rysunek{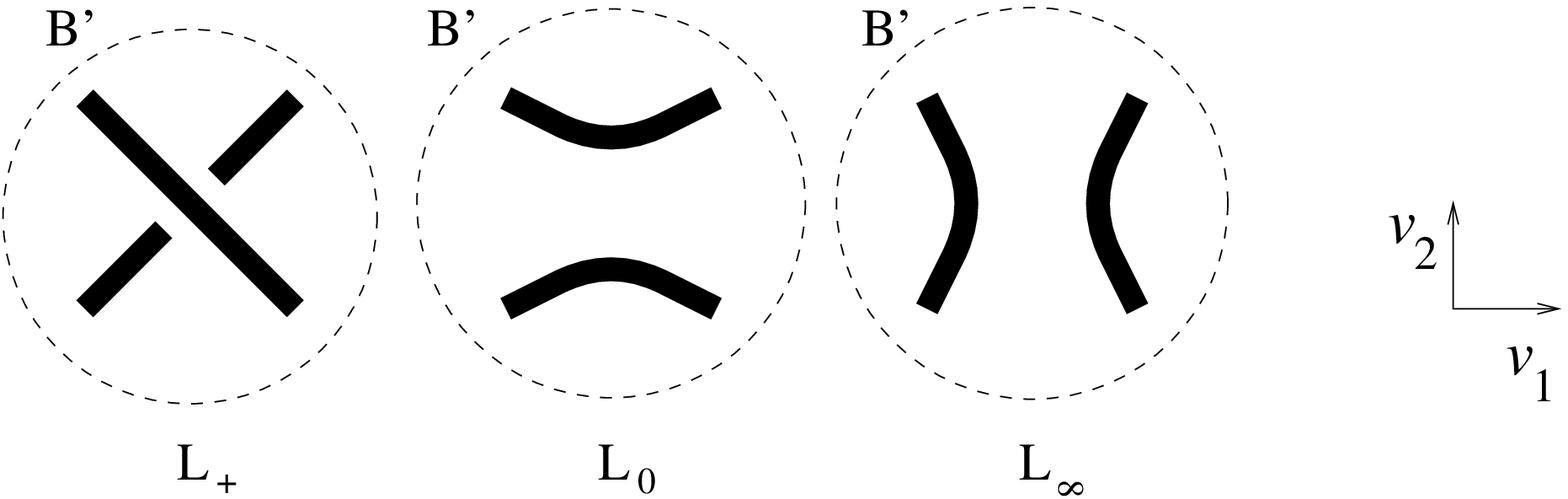}{3cm}}
\centerline{Figure 1}
\end{figure}

Notice that if a $3$-manifold $M$ is a disjoint union of manifolds
$M_1, M_2,...,M_n$ then \[{\cal S}_{2,\infty}(M;R,A)=
{\cal S}_{2,\infty}(M_1;R,A)\otimes_R {\cal S}_{2,\infty}(M_2;R,A)\otimes_R
... \otimes_R {\cal S}_{2,\infty}(M_n;R,A).\]
Therefore, we will assume, for simplicity, that all manifolds are connected
unless otherwise stated.

J. Barrett proved in [Bar] that the existence of a spin structure for any 
oriented $3$-manifold $M$ implies that for any knot $K$ in $M$ we can 
define its framing, $Spin(K)\in \Z_2,$ in such a way that the following 
theorem holds 

\begin{theorem}[\cite{Bar}]\label{2.2}\ \\
For any ring $R$ there exists an isomorphism of $R$-modules
$\phi:{\cal S}_{2,\infty}(M;R,A)\to {\cal S}_{2,\infty}(M;R,-A),$
such that $\phi([L])=(-1)^{\sum Spin(K)} [L],$ where the sum is over all
connected components of $L.$
\end{theorem}

The following fact, which is an easy generalization of results of 
[P-S-1], shows that skein modules ${\cal S}_{2,\infty}(M;R,A)$ are
particularly interesting for $A=1$ and $A=-1.$

\begin{fact}\label{2.3}\ \\
Let $M$ be any oriented manifold and let $R$ be a ring.
\begin{enumerate}

\item If $L_1$ and $L_2$ are two homotopic framed links or, equivalently,
if $L_2$ can be obtained from $L_1$ by an ambient isotopy and a sequence of
crossing changes, then $[L_1]=[L_2]$ in ${\cal S}_{2,\infty}(M;R,A)$ for
$A=1$ or $A=-1.$

\item $[L]\in {\cal S}_{2,\infty}(M;R,A),$ for $A=-1,$ does not depend on the
framing of $L.$

\item If $A=1$ or $A=-1$ then there is a multiplication operation on 
${\cal S}_{2,\infty}(M;R,A)$ uniquely determined by the following
condition:

For any two links $L_1,L_2$ in $M$ and any disjoint sum of them, 
$L_1\cup L_2,$ $[L_1]\cdot [L_2]=[L_1\cup L_2]$ in ${\cal S}_{2,\infty}
(M;R,-1).$ This multiplication is commutative, associative and has an 
identity $[\emptyset].$
\end{enumerate}
\end{fact}

We use Barrett's result to prove the following proposition.

\begin{proposition}\label{2.4}\ \\
The isomorphism of modules, $\phi,$ introduced in Theorem \ref{2.2}
is an isomorphism of $R$-algebras
$\phi:{\cal S}_{2,\infty}(M;R,-1)\to {\cal S}_{2,\infty}(M;R,1).$
\end{proposition}

\begin{proof}
Suppose that links $L,L'$ have connected components $K_1,K_2,...,K_n,$ and
$K_1',K_2',...,K_m'$ respectively. Then 
\[\phi([L][L'])=\phi([L\cup L'])=
(-1)^{\sum_{i=1}^n Spin(K_i)+\sum_{i=1}^m Spin(K_i')}[L\cup L']=\]
\[(-1)^{\sum_{i=1}^n Spin(K_i)}[L]\cdot (-1)^{\sum_{i=1}^n Spin(K_i)}[L']=
\phi([L])\phi([L']).\]
Therefore $\phi$ is a ring homomorphism and, by
Theorem \ref{2.2}, it is an isomorphism of $R$-algebras.
\end{proof}

Although the algebras ${\cal S}_{2,\infty}(M;R,-1)$ and ${\cal S}_{2,\infty}
(M;R,1)$ are isomorphic, the skein algebra ${\cal S}_{2,\infty}(M;R,-1)$ is 
much easier to understand (for example Fact \ref{2.3}(2) is not true for 
$A=1$).
Therefore, whenever we focus our attention on skein algebras we will prefer 
the choice $A=-1$ rather than $A=1.$ We believe that understanding
the properties of skein algebras ${\cal S}_{2,\infty}(M;R,-1)$ is the first 
step towards a better understanding of skein modules ${\cal S}_{2,\infty}
(M;R,A),$ with a generic $A\in R.$

We have proved in [P-S-1] that the skein module
${\cal S}_{2,\infty}(M;R,-1)$ depends only on the fundamental group of
the manifold $M.$ This motivated us to generalize the notion of skein
modules of manifolds to the notion of skein algebras of groups.

\begin{definition}\label{2.5}\ \\
Let us consider a group $G$ and a ring $R.$ 
Let ${\bf T}RG$ be the tensor algebra over the module $RG$.
Let ${\cal I}$ be an ideal of ${\bf T}RG$ generated by 
$e-2$ and expressions $g\otimes h -h\otimes g$, 
$g\otimes h - gh - gh^{-1}$, for any $g,h \in G.$
We define the skein algebra of $G$ as $
{\cal S}(G;R)={\bf T}RG/{\cal I}$.
\end{definition}

Obviously $\otimes$ is the multiplication in ${\cal S}(G; R)$
and usually $[gh]\ne [g]\otimes [h]$ in ${\cal S}(G; R).$
The algebra ${\cal S}(G; R)$ was denoted by ${\cal S}^-(G;R)$ in [P-S-1].

Very recently we have learned of work of K. Saito [S-1],[S-2].
He assigns for each group $G$ (and for a given ring $R$) {\it a universal 
character ring of $G.$} This ring is isomorphic to our skein algebra
${\cal S}(G;R).$

We list below a few elementary properties of ${\cal S}(G;R)$ which will 
be needed in the further part of the paper (See [P-S-1] for a proof).

\begin{fact}\label{2.6}\ \\
\begin{enumerate}
\item  For any $g\in G\  [g]=[g^{-1}]$ in ${\cal S}(G;R).$ 
\item  For any $g,h\in G\ [g]=[hgh^{-1}]$ and $[gh]=[hg]$ in ${\cal S}(G;R).$
\item  Any homomorphism of groups $f:G_1\to G_2$ yields a homomorphism
of R-algebras $f_*:{\cal S}(G_1;R)\to {\cal S}(G_2;R)$ such that
$f_*([g])=[f(g)]$ for any $g\in G.$
Moreover, epimorphisms of groups yield epimorphisms of algebras 
\item  (Universal Coefficient Property). Let $r: R \to R'$ be a
homomorphism of rings. We can think of $R'$ as an 
$R$-module. Then the identity map on $G$ induces an isomorphism of 
$R'$-algebras:
$${\overline r}: {\cal S}(G;R') \to {\cal S}(G;R)\otimes_R R'.$$
In particular, ${\cal S}(G;R) \simeq {\cal S}(G;\Z) \otimes_\Z R.$
\end{enumerate}
\end{fact}

We will also use another version of the Universal Coefficient Property
concerning skein modules ${\cal S}_{2,\infty}(M; R,A).$ The proof of
it is analogous to our proof of Fact \ref{2.6}(4) given in [P-S-1]
(Compare [P-1]).

\begin{fact}\label{2.7}\ \\
Let $r:R\to R'$ be a homomorphism of rings, $A^{\pm 1}\in R, A'=r(A).$
Then for any $3$-manifold $M$ there exists an isomorphism of $R'$-modules
$${\overline r}:{\cal S}_{2,\infty}(M;R',A') \to
{\cal S}_{2,\infty}(M;R,A)\otimes_R R'$$
such that ${\overline r}([L])=[L]$ for any framed link $L$ in $M.$
\end{fact}

We have proved the following theorem in [P-S-1] (see Lemmas 1.3 and 1.5).

\begin{theorem}\label{2.8}\ \\
For any 3-manifold $M$ and any ring $R$ there exists an isomorphism 
$$\hat{\xi}:{\cal S}_{2,\infty}(M;R,-1) \to {\cal S}(\pi_1(M);R),$$
such that for any knot $K$ in $M$ $\hat{\xi}([K])=-[\gamma],$ where 
$\gamma$ is an element of the fundamental group of $M$ represented by $K.$
\end{theorem}

\section{Group representations into $Sl_2(R)$}

In this section we are going to prove that under some general conditions
the skein algebra of a group $G$ is isomorphic to the algebra $TH_R(G)$ 
defined and investigated in the book by Brumfiel and Hilden, [B-H].

Let $G$ be any group and $R$ be a ring. Let $I$ be an ideal in the group ring
$RG$ generated by the expressions $g(h+h^{-1})-(h+h^{-1})g$ for $g,h\in G.$
We define (after Brumfiel and Hilden) an $R$-algebra $H_R(G)$ to be $RG/I.$

There exists an involution on $RG$ carrying $g$ to $g^{-1}$ for any $g\in G.$
Since this involution sends $I$ to $I,$ it yields an involution
$\imath : H_R(G)\rightarrow H_R(G)$ such that $\imath([g])=[g^{-1}].$

Let us consider the ring $M_2(R)$ of $2\times 2$ matrices (with coefficients
in $R$) and an involution $\imath:M_2(R)\rightarrow M_2(R)$ defined by
$$\imath\left(\matrix{a & b\cr c & d\cr}\right)=
\left(\matrix{d & -b\cr -c & a\cr}\right).$$

The following fact, observed by Brumfiel and Hilden (Proposition 1.3 [B-H]),
gives a good motivation for studying the algebras $H_R(G).$

\begin{fact}\label{3.1}\ \\
Any homomorphism $h:G\to Sl_2(R)$ extends to a homomorphism of $R$-algebras 
$\hat h: H_R(G)\rightarrow M_2(R)$ preserving involutions (i.e. 
$\imath\circ \hat h=\hat h\circ \imath$). Any involution preserving
homomorphism $\hat h: H_R(G)\rightarrow M_2(R)$ can be obtained in this way.
\end{fact}

Let $TH_R(G)$ be the subalgebra of $H_R(G)$ generated by elements of the form
$[g]+[g^{-1}], g\in G.$ One can easily notice (see [B-H]) that $TH_R(G)$ is a 
commutative algebra and $TH_R(G)$ is isomorphic, as an $R$-module, to 
$R+Span_R\{[g]+[g^{-1}]: g\in G\}\subset H_R(G).$

The following theorem shows that the skein algebra of a group $G$ is strongly
related to $TH_R(G).$

\begin{theorem}\label{3.2}\ \\
There exists an epimorphism of $R$-algebras, $\psi:{\cal S}(G;R)\rightarrow
TH_R(G),$ given by $\psi([g])=[g]+[g^{-1}],$ for any $g\in G.$
Moreover, if $2$ is not a zero divisor in ${\cal S}(G;R),$ then $\psi$ is
an isomorphism.\footnote{ The assumption that $2$ is not a zero divisor is
essential. For example, ${\cal S}(\Z_2;\Z_2)\simeq \Z_2[x]/(x^2)$ and
$TH_{\Z_2}(\Z_2)\simeq \Z_2.$}
\end{theorem}

\begin{proof}
Let $\psi_0:RG\rightarrow TH_R(G)$ be a homomorphism of $R$-modules given
by $\psi_0(g)=[g]+[g^{-1}]$ ($\psi_0$ is well defined because $RG$ is a free
$R$-module). We can uniquely extend $\psi_0$ to a homomorphism of $R$-algebras
$\psi_1$ defined on the tensor algebra of $RG,$
$\psi_1:TRG\rightarrow TH_R(G).$ Since $TH_R(G)$ is commutative,
$\psi_1(a\otimes b-b\otimes a)=0.$ Moreover,
$$\psi_1(e-2)=[e]+[e^{-1}]-2=0\quad {\rm and}$$
$$\psi_1(a\otimes b- ab- ab^{-1})=\psi_1(a)\psi_1(b)-\psi_1(ab)-
\psi_1(ab^{-1}) = ([a]+[a^{-1}])([b]+[b^{-1}])-$$
$$([ab]+[b^{-1}a^{-1}])-([ab^{-1}]+[ba^{-1}])=[a^{-1}]([b]+[b^{-1}])-
([b]+[b^{-1}])[a^{-1}]=0.$$ Therefore $\psi_1$ yields a homomorphism 
$\psi:{\cal S}(G;R)\rightarrow TH_R(G)$ such that $\psi([g])=[g]+[g^{-1}].$
Since $TH_R(G)$ is generated by the elements $[g]+[g^{-1}],$ for $g\in G,$
$\psi$ is an epimorphism.

We are going to prove that if $2$ is not a zero divisor in ${\cal S}(G;R),$
then $\psi$ is also a monomorphism. Before we start a proof of this
fact we will introduce two lemmas.

\begin{lemma}\label{3.3}\ \\
For any group $G$ and any ring $R$ the skein $R$-module ${\cal S}(G; R)$ is
generated by elements $[g],$ for $g\in G.$
\end{lemma}

\begin{proof}
By definition, ${\cal S}(G;R)$ is generated as an $R$-module by elements
$[g_1]\otimes [g_2]\otimes ...\otimes [g_n],$ where $g_1,g_2,...,g_n\in G.$
Moreover, $[g]\otimes [h]=[gh]+[gh^{-1}].$ Therefore any product
$[g_1]\otimes [g_2]\otimes ...\otimes [g_n]$ can be written as a linear 
combination of elements $[g_1\cdot g_2^{\pm 1}\cdot ... \cdot g_n^{\pm 1}].$
\end{proof}

Recall that $I\triangleleft RG$ was an ideal generated by expressions 
$g(h+h^{-1})-(h+h^{-1})g,$ for $g,h\in G.$ We will need the following lemma.

\begin{lemma}\label{3.4}\ \\
$I$ is, as an $R$-module, generated by expressions
$kg(h+h^{-1})-k(h+h^{-1})g,$ for $k,g,h\in~G.$
\end{lemma}

\begin{proof}
By definition, $I$ is spanned by the elements
\[k(g(h+h^{-1})-(h+h^{-1})g)l,\]
for $k,g,h,l\in G.$ But\\ 
$k(g(h+h^{-1})-(h+h^{-1})g)l=k[(g(h+h^{-1})l-gl(h+h^{-1})+gl(h+h^{-1})$\\
$-(h+h^{-1})gl]= kg((h+h^{-1})l-l(h+h^{-1}))+ k(gl(h+h^{-1})-(h+h^{-1})gl).$\\
\end{proof}

Let $\phi_0: RG\to {\cal S}(G;R)$ be a homomorphism of $R$-modules such that
$\phi_0(g)=[g]\in {\cal S}(G;R),$ for any $g\in G.$ Then, by Fact \ref{2.6}(2),
\[\phi_0(kg(h+h^{-1})-k(h+h^{-1})g)=[kgh]+[kgh^{-1}]-[khg]-[kh^{-1}g]=\]
\[[kgh]+[kgh^{-1}]-[gkh]-[gkh^{-1}]=[kg]\otimes [h]-[gk]\otimes [h]=0\] for
any $k,g,h\in G.$
Since $I$ is spanned by elements $kg(h+h^{-1})-k(h+h^{-1})g,$ we have
$\phi_0(I)=0.$
Therefore $\phi_0$ yields a homomorphism of $R$-modules
$$\phi: H_R(G)=RG/I \to {\cal S}(G;R).$$
Notice that $\phi\circ\psi([g])=\phi([g]+[g^{-1}])=[g]+[g^{-1}]=2[g]\in
{\cal S}(G;R).$ 
By Lemma 3.3 elements $[g], g\in G,$ generate ${\cal S}(G; R).$
Therefore $$\phi\circ\psi= 2\cdot id: {\cal S}(G;R)\to {\cal S}(G;R).$$
Hence if $2$ is not a zero divisor in ${\cal S}(G;R)$ then $\psi$ is a
monomorphism.
\end{proof}

Theorem \ref{3.2} implies the following corollary.

\begin{corollary}\label{3.5}\ \\
If one of the following conditions holds
\begin{enumerate}
\item [(1)] $\frac{1}{2}\in R;$ or
\item [(2)] ${\cal S}(G;R)$ is a free $R$-module and $2$ is not a 
zero divisor in $R$
\end{enumerate}
then $\psi:{\cal S}(G;R)\rightarrow TH_R(G)$ is an isomorphism.
\end{corollary}

The condition (2) of the above corollary motivated us to consider the 
following question: Under what conditions on $R$ and $G$ the skein
algebra $S(G,R)$ is a free $R$-module ?
There are only few classes of groups, $G,$ for which we are able
to prove that ${\cal S}(G;R)$ is free. Our proofs are based on 
topological methods.

\begin{theorem}\label{3.6}\ \\
If one of the following holds
\begin{enumerate}
\item [(1)] $G$ is an abelian group;
\item [(2)] $G$ is a free group;
\item [(3)] $G$ is the fundamental group of a surface;
\item [(4)] $G$ is the knot group of a $(2,2k+1)$-torus knot, $k\geq 0$;
\end{enumerate}
then the skein algebra ${\cal S}(G;R)$ is a free $R$-module for any ring $R.$
\end{theorem}

\begin{proof}
(1) We have proved in [P-S-1] (Lemma 2.2 and Theorem 2.3)
that if $G$ is abelian then ${\cal S}(G;\Z)$ is a free
$\Z$-module. Moreover, Fact \ref{2.6}(4) implies that ${\cal S}(G;R)=
{\cal S}(G;\Z)\otimes R.$ Therefore ${\cal S}(G;R)$ is a free $R$-module.\ \\
(3) We will see in Section 4 (Fact \ref{4.1}) that
the Kauffman bracket skein module of $F\times I,$ where $F$ is any 
surface, $I=[0,1],$ is free for any ring of coefficients $R$ and any 
$A^{\pm 1}\in R.$ Therefore ${\cal S}(\pi_1(F);R)\simeq
{\cal S}_{2,\infty}(F\times I; R,-1)$ is also free.\ \\
(2) It is a special case of (3).\ \\
(4) Let $M$ be the complement of a $(2,2k+1)$-torus knot in $S^3.$
D. Bullock proved in [B-5] that ${\cal S}_{2,\infty}(M; \Z[A^{\pm 1}],A)$ is 
a free $\Z[A^{\pm 1}]$-module. Let us consider any ring $R$ as a 
$\Z[A^{\pm 1}]$-algebra determined by a homomorphism 
$f:\Z[A^{\pm 1}]\rightarrow R, f(1)=1, f(A)=-1.$
By Fact \ref{2.7} and Theorem \ref{2.8},
${\cal S}(\pi_1(M);R)\simeq {\cal S}_{2,\infty}(M;R,-1)\simeq
{\cal S}_{2,\infty}(M;\Z[A^{\pm 1}],A)\otimes_{Z[A^{\pm 1}]} R$ is a free 
$R$-module.
\end{proof}

Theorems \ref{3.2}, \ref{3.6} and Corollary \ref{3.5} imply the 
following result.

\begin{corollary}\label{3.7}\ \\
If $G$ is abelian or if $G$ is a free group or a fundamental group of a 
surface or the group of a $(2,2k+1)$-torus knot then $TH_R(G)$ is a free 
$R$-module for any ring $R.$
\end{corollary}

If $G$ is a free group than one can give an alternative prove that 
$TH(G)$ is a free $R$-module by using Proposition 9.1[B-H]. For the case 
$G=F_3,F_4$ see also comments on page 9 and Proposition 6.6(iii), [B-H].

It is worth mentioning that the skein modules ${\cal S}_{2,\infty}(M; R,A)$ 
have much more complicated structure than the modules ${\cal S}(\pi_1(M); R).$
For example, if $M=S^1\times S^1\times S^1,$ then Theorems \ref{2.8} and
\ref{3.6} imply that ${\cal S}_{2,\infty}(M; R,-1)\simeq {\cal S}(\pi_1(M);R)=
{\cal S}(\Z^3;R)$ is a free $R$-module. We can prove, however, that the module 
${\cal S}_{2,\infty}(M; \Z[A^{\pm 1}],A)$ has a big torsion part and, in 
particular, is not free.

Another interesting example was considered by J. Hoste and the first author
in [H-P-3]. They proved that if $M$ is a classical Whitehead manifold then 
the skein module ${\cal S}_{2,\infty}(M; \Z[A^{\pm 1}],A)$ is torsion free 
and infinitely generated. On the other hand Theorem \ref{2.8} and the 
definition of skein algebra imply that ${\cal S}_{2,\infty}(M; R,-1)\simeq 
{\cal S}(\{e\}; R)\simeq R$ for any ring $R.$

\section {On skein algebras of surfaces}

Let $F$ be an oriented surface and $I$ be the closed interval $[0,1].$
In this section we are investigating skein modules
${\cal S}_{2,\infty}(F\times I; R,A).$ They admit a natural operation of 
multiplication which can be defined in the following way.

For any two framed links $L_1,L_2$ in $F\times I$ we define
$L_1\cdot L_2$ to be a union of $L_1$ and $L_2,$ where $L_1$ is isotopically 
pushed into $F\times (1/2,1)$ and $L_2$ is isotopically pushed into 
$F\times (0,1/2).$ This multiplication is well defined
on the set of ambient isotopy classes of links in
$F\times I$ because if links $L_i$ and $L_i'$ are ambient
isotopic in $F\times I,$ for $i=1,2,$ then $L_1\cdot L_2$ is ambient
isotopic to $L_1'\cdot L_2'.$ To show this, note that if links $L$ and
$L'$ are ambient isotopic in $F\times I$ then we can choose an ambient
isotopy between them which fixes the boundary of $F\times I.$
Therefore, we can obtain an ambient isotopy between $L_1\cdot L_2$ and
$L_1'\cdot L_2'$ by composing an ambient isotopy between $L_1$ and $L_1'$ in
$F\times [1/2,1],$ which fixes $\partial(F\times [1/2,1]),$ with an
an ambient isotopy between $L_2$ and $L_2'$ in
$F\times [0,1/2],$ which fixes $\partial(F\times [0,1/2]).$

Having defined a product on the set of all ambient isotopy classes
of framed, unoriented  links in $F\times I,\ {\cal L}_{fr}(F\times I),$
we can linearly extend it to a multiplication in the module 
$R{\cal L}_{fr}(F\times I)$ and, finally, obtain a multiplication
in ${\cal S}_{2,\infty}(F\times I; R,A).$ The element $[\emptyset]=1$ is the 
identity in ${\cal S}_{2,\infty}(F\times I; R,A).$
If $A=-1$ then the multiplication defined above coincides with the 
multiplication considered in Sections 2 and 3.

Notice that the multiplication in ${\cal S}_{2,\infty}(F\times I; R,A)$
is uniquely determined by the surface $F$ but it is not necessarily 
unique to the manifold $F\times I.$ To illustrate this phenomenon let
us consider a disc with two holes, $F_0,$ and a punctured torus, $F_1.$
Using Fact \ref{4.1} one can easily prove that 
${\cal S}_{2,\infty}(F_1\times I;R,A)$ is a non-commutative algebra.
On the other hand ${\cal S}_{2,\infty}(F_0\times I; R,A)$ is
commutative by Corollary \ref{4.4}. Hence the algebras
${\cal S}_{2,\infty}(F_1\times I;R,A)$ and 
${\cal S}_{2,\infty}(F_0\times I; R,A)$ are clearly not isomorphic
even though there is a homeomorphism of topological spaces
$F_0\times I\simeq F_1\times I.$

In the first part of this section we prove a few basic properties of skein 
algebras of surfaces. We show that ${\cal S}_{2,\infty}(F\times I; R,A)$
can be considered as an algebra over a ring of polynomials in many variables.
We also describe ${\cal S}_{2,\infty}(F\times I; R,A)$ for $F$ being a 
$d$-punctured sphere, where $d=0,1,2,3.$
Moreover, we announce two more sophisticated results on skein algebras
${\cal S}_{2,\infty}(F\times I; R,A).$
We claim that if $R$ has no zero divisors then the algebra
${\cal S}_{2,\infty}(F\times I; R,A)$ has no zero divisors either.
Moreover, ${\cal S}_{2,\infty}(F\times I; R,A)$ considered as an algebra over 
a ring of polynomials is a central algebra. These two results require
relatively long proofs based on Dehn's theorem classifying curves on surfaces.
For that reason we decided to publish proofs of these results in a separate 
paper, [P-S-2].

In this section we assume that $F$ is an oriented, compact, connected
2-manifold. We identify $F$ with $F\times \{1/2\}\subset F\times I.$
Let $L(F)$  be the set of all unoriented links in $F$ without (homotopically) 
trivial components. We assume, as usual, that $\emptyset\in L(F).$ Notice that
any link in $L(F)$ has a natural framing i.e. the framing parallel to $F.$
Therefore we can consider elements of $L(F)$ as framed links in $F\times I.$
Let ${\cal B}(F)$ be a set containing exactly one link from each class of 
ambient isotopic links in $L(F).$
We start with the following basic result (See [P-1], [H-P-2]).

\begin{fact}\label{4.1}\ \\
For any surface $F,$ any ring $R,$ and any element $A^{\pm 1}\in R$ the skein 
module ${\cal S}_{2,\infty}(F\times I; R,A)$ is a free $R$-module with a basis
$\{[L]\}_{L\in {\cal B}(F)}.$
\end{fact}

For any surface $F\not \simeq S^1\times I$ we define $d(F)$ to be the number
of boundary components of $F.$ For $F=S^1\times I$ we set $d(F)=1.$

\begin{lemma}\label{4.2}\ \\
Let $F\not \simeq S^1\times I$ be a surface with $d$ boundary
components denoted by $\partial_1F,\partial_2F,...,\partial_dF.$ 
Let $N_i,$ for $i\in \{1,2,...,d\},$ be a regular neighborhood
of $\partial_iF.$ We assume that $N_1,N_2,...,N_d$ are disjoint.
Let $K_i$ be a knot in $N_i(F)\subset F$ parallel to $\partial_iF.$
For $F \simeq S^1\times I$ set $d=1,\ N_1=F$ and define $K_1$ to
be a knot in $F$ parallel to $\partial F.$
Then
\begin{enumerate}
\item There is a natural isomorphism of $R$-algebras
$$f:R[\partial_1,\partial_2,...,\partial_d]\to {\cal S}_{2,\infty}
(N_1\times I\cup N_2\times I\cup ...\cup N_d\times I; R, A)$$ such that
$f(\partial_i)=[K_i].$

\item The embedding $i:N_1\times I\cup N_2\times I\cup ... \cup N_d\times I 
\to F\times I$ induces a monomorphism $i_*:R[\partial_1,\partial_2,...,
\partial_d]\to {\cal S}_{2,\infty}(F\times I; R,A).$
\end{enumerate}
\end{lemma}

\begin{proof}
\begin{enumerate}
\item Notice that there is an isomorphism between
${\cal S}_{2,\infty}(N_1\times I\cup N_2\times I\cup ...
\cup N_d\times I; R, A)$ and ${\cal S}_{2,\infty}(N_1\times I;R,A)\otimes
{\cal S}_{2,\infty}(N_2\times I;R,A)\otimes ... \otimes
{\cal S}_{2,\infty}(N_d\times I;R,A)$ matching elements of the form
\[[L_1\cup L_2\cup ...\cup L_d]\in {\cal S}_{2,\infty}(N_1\times I\cup 
N_2\times I\cup ...\cup N_d\times I; R, A),\] 
where $L_i$ is a link in $N_i, i\in \{1,2,...,d\},$ with elements 
\[[L_1]\otimes [L_2]\otimes ...
\otimes [L_d]\in {\cal S}_{2,\infty}(N_1\times I;R,A)\otimes
... \otimes {\cal S}_{2,\infty}(N_d\times I; R, A).\]
Since $N_i\simeq S^1\times I,$ we can assume that ${\cal B}(N_i)=
\{K_i,K_i^2,K_i^3,...\},$ where $K_i^n$ denotes a link in $N_i$ composed 
of $n$ parallel copies of $K_i.$ By Fact \ref{4.1} there is an isomorphism of 
$R$-modules $f_i: R[\partial_i]\to {\cal S}_{2,\infty}(N_i\times I)$ such that
$f_i(\partial_i^n)=[K_i^n].$ One can easily see that $f_i$ is also an 
isomorphism of $R$-algebras. Therefore
\[f=f_1\otimes f_2\otimes ...\otimes f_n: R[\partial_1,\partial_2,...,
\partial_d]= R[\partial_1]\otimes R[\partial_2]\otimes ...\otimes 
R[\partial_d]\quad \to\]
\[{\cal S}_{2,\infty}(N_1\times I;R,A)\otimes
{\cal S}_{2,\infty}(N_2\times I;R,A)\otimes ... \otimes
{\cal S}_{2,\infty}(N_d\times I;R,A)\simeq\]
\[{\cal S}_{2,\infty}(N_1\times I\cup N_2\times I
\cup ... \cup N_d\times I; R, A)\] is an isomorphism.

\item We have defined $K_1,...,K_d$ in such a way that any two links
$K_1^{n_1}\cup K_2^{n_2}\cup ...\cup K_d^{n_d}$ and 
$K_1^{n_1'}\cup K_2^{n_2'}\cup ...\cup K_d^{n_d'}$ are ambient isotopic in 
$F\times I$ if and only if $(n_1,n_2,...,n_d)=(n_1',n_2',...,n_d').$ 
Therefore we can assume that all links of the form $K_1^{n_1}\cup K_2^{n_2}
\cup ...\cup K_d^{n_d},$ for any $n_1,n_2,...,n_d\in \{0,1,2,...\}$ belong to 
the set ${\cal B}(F).$

Let us consider a basis of $R[\partial_1,\partial_2,...,\partial_d]$ composed
of all monomials. 
Notice that $i_*(\partial_1^{n_1}\partial_2^{n_2}...
\partial_d^{n_d})=[K_1^{n_1}\cup K_2^{n_2}\cup ...\cup K_d^{n_d}]\in
{\cal S}_{2,\infty}(F\times I;R,A).$
Therefore $i_*$ carries different elements of the basis of $R[\partial_1,
\partial_2,...,\partial_d]$ to different elements of the basis ${\cal B}(F)$ 
of the module ${\cal S}_{2,\infty}(F\times I;R,A).$ Hence $i_*$ is a 
monomorphism.
\end{enumerate}
\end{proof}

Let ${\cal B'}(F)$ be a set of all links in ${\cal B}(F)$ without components 
parallel to the boundary of $F.$ Notice that $\emptyset\in {\cal B'}(F).$ 
Fact \ref{4.1} and Lemma \ref{4.2} imply the following corollary.

\begin{corollary}\label{4.3}\ \\
${\cal S}_{2,\infty}(F\times I;R,A)$ considered as an 
$R[\partial_1,\partial_2,...,\partial_{d(F)}]$-algebra has a basis 
$\{[L]\}_{L\in {\cal B'}(F)}.$
\end{corollary}

The next corollary describes the skein algebra
${\cal S}_{2,\infty}(F\times I;R,A)$ for a few important surfaces $F$
(Compare [B-P]).

\begin{corollary}\label{4.4}\ \\
Let $F$ be a sphere or a disc or an annulus or a disc with two 
holes. Then $i_*: R[\partial_1,...,\partial_{d(F)}]\to {\cal S}_{2,\infty}
(F\times I;R,A)$ is an isomorphism of $R$-algebras.
\end{corollary}

\begin{proof}
Notice that ${\cal B'}(F)=\{\emptyset\}$ for surfaces $F$ listed in Corollary
\ref{4.4}. Therefore the statement of Corollary \ref{4.4} is implied by 
Corollary \ref{4.3}.
\end{proof}

If we combine Corollary \ref{4.4} with Theorem \ref{2.8} we will get the 
following result.

\begin{corollary}\label{4.5}\ \\
The following homomorphisms
\begin{enumerate}
\item $f_1:R[\partial_1]\to {\cal S}(\Z;R), f_1(\partial_1)=[1]$
\item $f_2:R[\partial_1,\partial_2,\partial_3]\to {\cal S}(F_2; R), 
f_2(\partial_1)=[a], f_2(\partial_2)=[ab], f_2(\partial_3)=[b]$
\end{enumerate}
are isomorphisms of $R$-algebras.
\end{corollary}

Proofs of the following two crucial results about skein algebras of surfaces
will appear in our paper [P-S-2].

\begin{theorem}\label{4.6}\ \\
Let $A\in R$ be an invertible element which is not a root of unity in $R.$
Then ${\cal S}_{2,\infty}(F\times I;R,A)$ is a central
$R[\partial_1,...,\partial_{d(F)}]$-algebra i.e.
$i_*(R[\partial_1,...,\partial_{d(F)}])$ is the center of the algebra
${\cal S}_{2,\infty}(F\times I;R,A).$
\end{theorem}

\begin{theorem}\label{4.7}\ \\
Let $R$ be a ring without zero divisors. 
Suppose that
\begin{enumerate}
\item $F$ is an orientable surface, $M=F\times I,$ and $A^{\pm 1}\in R;$ or
\item $F$ is an unorientable surface of an even, negative Euler 
characteristic, $M$ is a twisted $I$ bundle over $F,$
and $A=\pm 1.$
\end{enumerate}
Then the skein algebra ${\cal S}_{2,\infty}(F\times I; R,A)$ has no zero 
divisors.
\end{theorem}

The above theorem is not true for a Klein bottle, $KB.$ One can check that
if ${1\over 2}\in R$ then ${\cal S}_{2,\infty}(KB\times I; R,-1)\simeq 
{\cal S}(<a,b\ |\ aba^{-1}b>,R)\simeq R[x,y]/(x(y^2-4))$ has zero divisors.

We will observe in Section 7 that Theorem \ref{4.7} has important 
applications to the theory of $Sl_2(\C)$-character varieties.

\section{Relative skein algebras}

In this section we will define relative skein modules of
$3$-manifolds. The relative skein module of $F\times I,$ where $F$ is
an oriented surface, admits an algebra structure and therefore it is
called a relative skein algebra of $F.$
We will prove that the Brumfiel and Hilden algebra $H(G),$ for
$G=\pi_1(F),$ is isomorphic to a relative skein algebra of $F.$

Let $M$ be an oriented $3$-manifold with a boundary
and let $\gamma_0,\gamma_1:[0,1]\to \partial M$ be two fixed simple arcs in
the boundary of $M.$ By a special framed arc in $M$ we mean an embedding
$\imath: [0,1]\times I\to M$ such that
$\imath(s,0)=\gamma_0(s)\in \partial M, \imath(s,1)=\gamma_1(s)\in
\partial M,$ for any $s\in [0,1].$ We also assume that
$\imath(s,t)$ lies in the interior of $M$ for $t\in (0,1).$
We can consider a special framed arc as a ribbon inside of $M$ whose
ends lie exactly on $\gamma_0([0,1])$ and $\gamma_1([0,1]).$

We define a relative framed link in $M$ to be a disjoint union of a 
special framed arc with a framed unoriented link lying in the interior
of $M.$ Since we consider $\emptyset$ as a link, any special framed
arc in $M$ is also a relative framed link. We say that two relative
framed links $L$ and $L'$ are ambient isotopic if there is an ambient
isotopy of $M$ which carries $L$ to $L'$ and is fixed on 
$\partial M.$

Let ${\cal L}^{rel}_{fr}(M)$ denote the set of all ambient isotopy 
classes of relative framed links in $M.$ 
We define a relative skein module of $M$ in the same way as
it was done in Definition \ref{2.1} but we replace ${\cal L}_{fr}(M)$ by 
${\cal L}_{fr}^{rel}(M).$

\begin{definition}\label{5.1}\ \\
Let $M$ be a $3$-manifold and let $R$ be any ring with a specified 
invertible element $A.$ Let $S_{2,\infty}^{rel}$ be the submodule of 
$R{\cal L}_{fr}^{rel}(M)$ generated by 
\begin{enumerate}
\item skein expressions $L_+-AL_0 - A^{-1}L_{\infty}$ where 
$L_+, L_0, L_{\infty}$ are any skein related framed relative links in 
$M$ \footnote{Skein related links were introduced in Definition
\ref{2.1}.}.

\item $L\cup \bigcirc +(A^2+A^{-2})\cdot L$ for any relative framed link
$L$ in $M.$
\end{enumerate}

We define the relative Kauffman bracket skein module of $M$ as the quotient 
$${\cal S}_{2,\infty}^{rel}(M;R,A)= R{\cal L}_{fr}^{rel}(M)/
S_{2,\infty}^{rel}.$$
\end{definition}

The definition of a relative framed link in $M,$ as well as
the definition of a relative skein module of $M,$ depends on
a particular choice of simple arcs $\gamma_0,\gamma_1: [0,1]\to
\partial M.$ One can easily see, however, that different choices of 
arcs $\gamma_0,\gamma_1$ give isomorphic relative skein modules.

In this section we will be interested in the case when $M=F\times I,$
for a surface $F.$ We will notice that there is a
natural product on ${\cal S}_{2,\infty}^{rel}(F\times I;R,A)$ similar to that
on ${\cal S}_{2,\infty}(F\times I;R,A).$
However the relative skein modules of cylinders over surfaces,
${\cal S}_{2,\infty}^{rel}(F\times I ;R,A),$
are usually more complicated than the classical skein modules,
${\cal S}_{2,\infty}(F\times I;R,A).$ We are able to describe them only for 
surfaces with boundaries. Our description is analogous to that given in 
Fact \ref{4.1}

As before we identify $F$ with $F\times \{1/2\}\subset F\times I.$
Let us assume that $\gamma_0,\gamma_1:[0,1]\to \partial F\subset
\partial(F\times I)$ are two disjoint arcs lying in the same component of
$\partial F.$ Let $p_0=\gamma_0(1/2),p_1=\gamma_1(1/2).$
By a relative (unframed) link in $F$ we will mean a disjoint union
$L=L_0\cup K,$ where $L_0$ is a link in the interior of $F$ and
$K:[0,1]\to F$ is an arc, $K(0)=p_0,\ K(1)=p_1,$ and $K(t)$ lies in
the interior of $F,$ for $t\in (0,1).$
We denote the set of all relative links in $F$ without (homotopically) trivial
components by $L^{rel}(F).$
We assume that $\emptyset\in L^{rel}(F).$  Let ${\cal B}^{rel}(F)$ be a set
containing exactly one link from each class of ambient isotopic relative links
in $L^{rel}(F).$

Notice that a regular neighborhood in $F$ of any relative link 
$L=L_0\cup K\in L^{rel}(F)$ is a framed relative link in $F\times I.$
Therefore one can identify elements of ${\cal B}^{rel}(F)$ with framed relative
links in ${\cal S}_{2,\infty}^{rel}(F\times I;R,A).$
 
\begin{theorem}\label{5.2}\ \\
Let $F$ be a surface with a boundary. Then
${\cal S}_{2,\infty}^{rel}(F\times I;R,A)$ is a free $R$-module with a basis
$\{[L]\}_{L\in {\cal B}^{rel}(F)}.$
\end{theorem}

{\it Sketch of proof:}
Let $L$ be a relative framed link in $F\times I.$ Then $L$ can be represented 
by a diagram $D$ in $F.$ Notice that if we smooth each of the crossings of $D$
in two possible ways, we will obtain a collection of relative links in $F.$ 
This implies that $L$ can be represented by a linear combination of links 
$\sum r_i[L_i],$ where $L_i\in {\cal B}^{rel}(F),\ r_i\in R.$

Two diagrams, $D_1, D_2\subset F$ represent the same relative framed link
$L\subset F\times I$ if one can be transformed into the other by a sequence of
Reidemeister moves. One can check that these 
moves do not change the linear combination $\sum r_i[L_i]$ assigned to $L.$
Therefore, there is a homomorphism of $R$-modules
$$f_0: R{\cal L}_{fr}^{rel}(F\times I)\to R{\cal B}^{rel}(F),$$ such that
$f_0(L)=\sum r_i[L_i].$ One can show that this homomorphism yields an
isomorphism $f_1: {\cal S}_{2,\infty}^{rel}(F\times I;R,A)\to
R{\cal B}^{rel}(F).$
\Box

\medskip
Theorem \ref{5.2} is not true for closed surfaces. 
Let $F=S^2$ and let $L$ be any framed arc in $F\times I$ whose ends lie on 
$\gamma([0,1])\times \{0\}$ and $\gamma([0,1])\times \{1\},$ for some curve 
$\gamma:[0,1]\to F.$ Then one can show that $[L]\ne 0$ in ${\cal S}_
{2,\infty}^{rel}(F\times I;R,A)$ but $(A^6-1)[L]=0.$ Therefore 
${\cal S}_{2,\infty}^{rel}(S^2\times I;R,A)$ is not free.

In order to define a multiplication operation on
${\cal S}_{2,\infty}^{rel}(F\times I;R,A)$ it is convenient to
consider arcs $\gamma_0,\gamma_1$ defined as follows.
Let $\gamma:[0,1]\to F$ be a simple arc in the interior of $F$ and let
$\gamma_0,\gamma_1\in F\times I$ be given by
$\gamma_0(t)=(\gamma(t),0), \gamma_1(t)=(\gamma(t),1),$ for any $t\in [0,1].$
Let $h_1: F\times [0,1] \to F\times [{1 \over 2},1], h_2: F\times [0,1]\to
F\times [0,{1 \over 2}]$ be maps given by $h_1(x,t)=(x, {1\over 2}+
{1\over 2}t), h_2(x,t)=(x,{1\over 2}t)$ for $x\in F, t\in [0,1].$ We define
the product of relative framed links $L_1, L_2\subset F\times I,$ to be
$L_1\cdot L_2 = h_1(L_1)\cup h_2(L_2)\subset F\times [{1\over 2},1]
\cup F\times [0,{1\over 2}]= F\times [0,1].$ One can easily see that
$h_1(L_1)\cup h_2(L_2)$ is a relative framed link in $F\times I.$
Moreover, we can uniquely extend the above multiplication to a multiplication
operation on ${\cal S}^{rel}_{2,\infty}(F\times I; R,A).$ This
multiplication will be associative but (in general) not commutative.
The module ${\cal S}^{rel}_{2,\infty}(F\times I; R,A)$ considered
together with this multiplication will be called the
relative skein algebra of a surface $F.$ Whenever we will talk about
the relative skein algebras of surfaces we will always assume that the
arcs $\gamma_0,\gamma_1$ are defined as above. Notice that the
any choice of $\gamma:[0,1]\to F$ gives the same (up to an
isomorphism) algebra structure on ${\cal S}^{rel}_{2,\infty}(F\times I; R,A).$

The next result shows that the relative skein algebra of $F$ can be
considered as a generalization (or a deformation) of
the algebra $H_R(G),\ G=\pi_1(F),$ introduced in Section 3. 

\begin{theorem}\label{5.3}\ \\
Let $F$ be any oriented surface and let $R$ be a ring.
Then $H_R(\pi_1(F))$ is isomorphic as an $R$-algebra to
${\cal S}^{rel}_{2,\infty}(F\times I; R,-1).$
\end{theorem}

\begin{proof}
Let us fix a point $p\in F$ and assume that $F=F\times \{1/2\}\subset 
F\times I,\ I=[0,1].$ Let $G=\pi_1(F,p).$
We will define a homomorphism $f_0:RG\to {\cal S}^{rel}_{2,\infty}
(F\times I;R,-1)$ so that $f_0$ vanishes on
$I=(h(g+g^{-1})-(g+g^{-1})h)\triangleleft RG,$ and thus factors
through $H_R(G).$

In the proof we will consider relative (unframed) links in $F\times I.$
A relative (unframed) link in $F\times I$ is a disjoint union of 
a (possibly empty) unoriented link in $F\times I$ with an arc in 
$F\times I$ joining $(p,0)$ with $(p,1).$
We denote the set of all ambient isotopy classes of relative links in
$F\times I$ by ${\cal L}^{rel}(F\times I).$

Notice that $[L]\in {\cal S}^{rel}_{2,\infty}(F\times I; R,-1)$ does not depend on
the framing of a relative link $L.$ Therefore we can define
${\cal S}^{rel}_{2,\infty}(F\times I; R,-1)$ using unframed links,
${\cal S}^{rel}_{2,\infty}(F\times I; R,-1)= R{\cal L}^{rel}(F\times I)/
{\cal S}^{rel},$ where ${\cal S}^{rel}$ is a submodule of
$R{\cal L}^{rel}(F\times I)$ generated by expressions
\begin{enumerate}
\item $L_+ +L_0+L_{\infty},$ for any skein related (unframed) relative links
$L_+,L_0,L_{\infty}$ in $F\times I.$
\item $L\cup \bigcirc +2\cdot L,$ for any (unframed) relative link in
$L$ in $F\times I.$
\end{enumerate}

Suppose that $g\in G=\pi_1(F,p)$ is the homotopy class of an oriented loop
$\xi: [0,1]\to F, \xi(0)=\xi(1)=p.$ Then we define an arc
$K_g:[0,1]\to F\times I, K_g(t)=(\xi(t),1-t),\ t\in [0,1].$
Notice that $K_e$ is the trivial arc, $K_e: [0,1]\to F\times I, 
K_e(t)=(p,1-t).$

Let $f_0: RG\to {\cal S}^{rel}_{2,\infty}(F\times I; R,-1)$ be an $R$-module
homomorphism such that $f_0(g)=[K_g]\in {\cal
S}^{rel}_{2,\infty}(F\times I; R,-1).$ 
Notice that $f_0(g_1g_2)=[K_{g_1g_2}]=
[K_{g_1}][K_{g_2}]=f_0(g_1)f_0(g_2).$ Therefore $f_0$ is a 
homomorphism of $R$-algebras.

Let $L_+=K_e\cup K,$ where $K_e$ is the trivial arc and $K$ is a knot in 
$F\times I$ representing (with a given orientation) the conjugacy 
class of some $g\in G.$ Homotopic relative framed links 
in $F\times I$ are identified in
${\cal S}^{rel}_{2,\infty}(F\times I; R,-1).$ Therefore we can
place $K$ as in Fig. 2. Then $[L_+]=-[L_0]-[L_{\infty}]$ in
${\cal S}^{rel}_{2,\infty}(F\times I; R,-1),$ where $L_0=K_g,
L_{\infty}=K_{g^{-1}}$ are presented in Figure 2.

\begin{figure}
\centerline{\rysunek{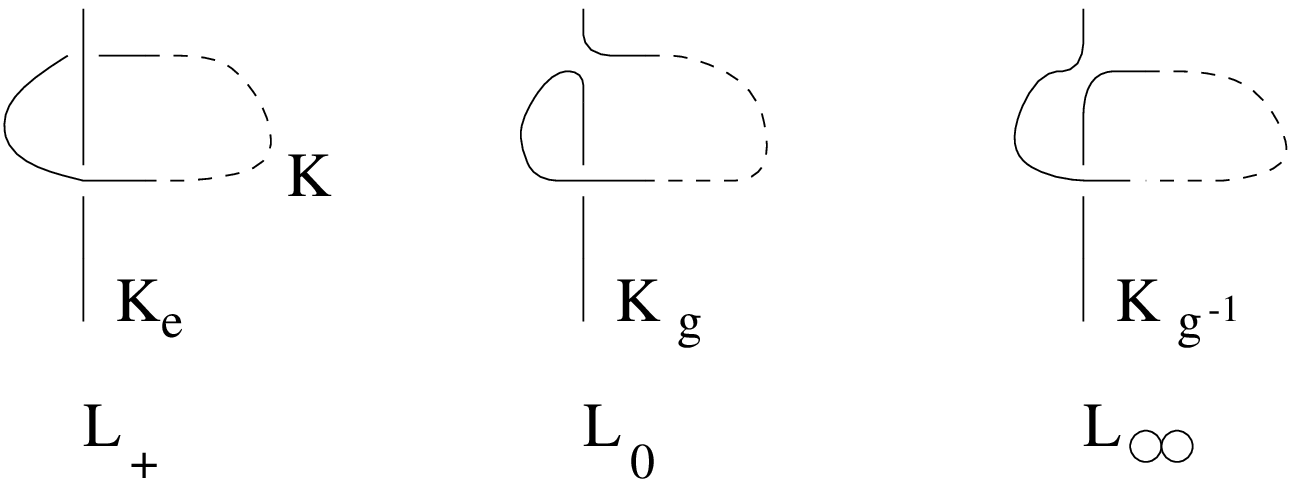}{3.2cm}}
\centerline{Figure 2}
\end{figure}

Hence $[L_+]=-f_0(g+g^{-1}).$ Notice that $[L_+]$ commutes with any relative 
link in $F\times [0,1].$ Therefore $f_0(h(g+g^{-1}))=f_0((g+g^{-1})h)$ for any
$g,h\in G$ and $f_0$ vanishes on $I=(h(g+g^{-1})-(g+g^{-1})h)\triangleleft RG.$
Hence $f_0$ induces a homomorphism $f_1:H_R(G)\to 
{\cal S}^{rel}_{2,\infty}(F\times I; R,-1).$

\begin{lemma}\label{5.4}\ \\
$f_1$ is an epimorphism.
\end{lemma}

\begin{proof}
Notice that for any relative link $L$ in $F\times I,$ which is
composed of an arc $K_g$ and knots $L_1,L_2,...,L_n,\ L=K_g\cup L_1\cup
L_2\cup ... \cup L_n,$ we have $[L]=[K_g][L_1\cup K_e]\cdot
[L_2\cup K_e]\cdot ... \cdot [L_n\cup K_e]\in {\cal S}^{rel}_{2,\infty}
(F\times I;R,-1).$
Let $L_i$ with a given orientation correspond to the conjugacy class
of some $g_i\in G.$ Then, as we have seen before,
$[L_i\cup K_e]= -f_1(g_i+g_i^{-1}).$ Hence $[L_i\cup K_e]\in Im\ f_1.$
Since $[K_g]=f_1(g)$ and $f_1$ is a homomorphism of algebras, 
$[L]\in Im\ f_1.$ The elements $[L],$ where $L$ are relative links,
generate ${\cal S}^{rel}_{2,\infty}(F\times I; R,-1).$ Therefore $f_1$ is an
epimorphism.
\end{proof}

We are going to complete the proof of Theorem \ref{5.3} by showing that 
$f_1$ is a monomorphism.

Let $J\triangleleft RG\otimes_R {\cal S}(G;R)$ be an ideal generated by 
elements $g\otimes [h]- gh\otimes 1 - gh^{-1}\otimes 1,$ for any $g,h\in G.$ 
Let $\rho_0: R{\cal L}^{rel}(F\times I)\to (RG\otimes_R {\cal S}(G;R))/J$ be 
a homomorphism of $R$-modules defined in the following way:

If $L$ is a relative link in $F\times I$ of the form $L=K_g\cup L',$ where 
$L'$ is a link disjoint from an arc $K_g,$ then
$\rho_0(L)=g\otimes \hat{\xi}([L']) + J$ \footnote{The homomorphism 
$\hat{\xi}$ was introduced in Theorem \ref{2.8}.}.

\begin{lemma}\label{5.5}\ \\
Let $L_+,L_0,L_{\infty}$ be any skein related links in $F\times I.$
Then $\rho_0(L_++L_0 +L_{\infty})=0.$
\end{lemma}

\begin{proof}
Let $L_+=K_g\cup L',$ where $g\in G$ and $L'$ is a link in $F\times I.$
We will call the crossing in $L_+$ which has to be smoothed in order to
obtain $L_0$ or $L_{\infty}$ {\it a specified crossing} in $L_+.$
There are three possibilities:
\begin{enumerate}

\item The specified crossing in $L_+$ 
is between two connected components of $L'$ or is a self-crossing of a
component of $L'.$

Then $L_0=K_g\cup L_0', L_{\infty}=K_g\cup L_{\infty}'$ and
$\rho_0(L_++L_0+L_{\infty})=g\otimes \hat{\xi}([L'])+
g\otimes \hat{\xi}([L_0'])+g\otimes \hat{\xi}([L_{\infty}'])+J=
g\otimes \hat{\xi}([L']+[L_0']+[L_{\infty}'])+J.$
Since $L',L_0',L_{\infty}'$ are skein related, $[L']+[L_0']+[L_{\infty}]=0$ in
${\cal S}_{2,\infty}(F\times I; R,-1).$ Therefore
$\rho_0(L_++L_0+L_{\infty})=0.$

\item The specified crossing in $L_+$ is a self-intersection of
$K_g.$

Then $L_0=K_{h_1}\cup L'\cup L'',\ L_{\infty}=K_{h_1h_2^{-1}}\cup L'$ have
forms presented in Fig. 3, $h_1h_2=g,$ and the knot $L''$ corresponds (with 
a given orientation) to the conjugacy class of $h_2\in G.$ 

\begin{figure}
\centerline{\rysunek{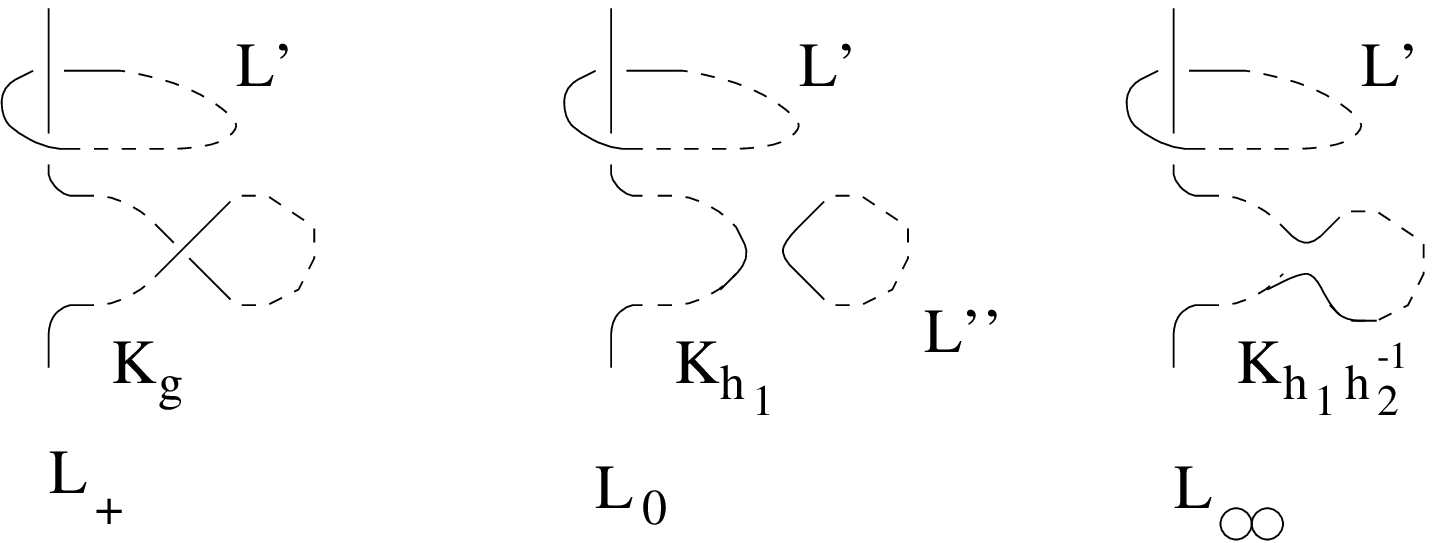}{3.5cm}}
\centerline{Figure 3}
\end{figure}

Therefore
$\rho_0(L_++L_0+L_{\infty})=g\otimes \hat{\xi}([L'])+h_1\otimes \hat{\xi}
([L'\cup L'']) + h_1h_2^{-1}\otimes \hat{\xi}([L'])+J=
(g\otimes 1+h_1\otimes (-[h_2])+h_1h_2^{-1}\otimes 1+ J)\cdot (1\otimes 
\hat{\xi}([L']) +J)=(g\otimes 1 -h_1h_2\otimes 1 -h_1h_2^{-1}\otimes 1+
h_1h_2^{-1}\otimes 1 +J)\cdot (1\otimes \hat{\xi}([L'])+J)= 0.$

\item The specified crossing in $L_+$ is between $K_g$ and a connected
component $K$ of $L'.$

Choose some orientation of $K.$ Then $K$ corresponds to the conjugacy 
class of some $h\in \pi_1(F,p)=G.$ Let $L''=L'\setminus K.$ Then $L_0$ and 
$L_{\infty}$ have forms $K_{gh}\cup L''$ and $K_{gh^{-1}}\cup L''$ 
as shown in Figure 4.

\begin{figure}
\centerline{\rysunek{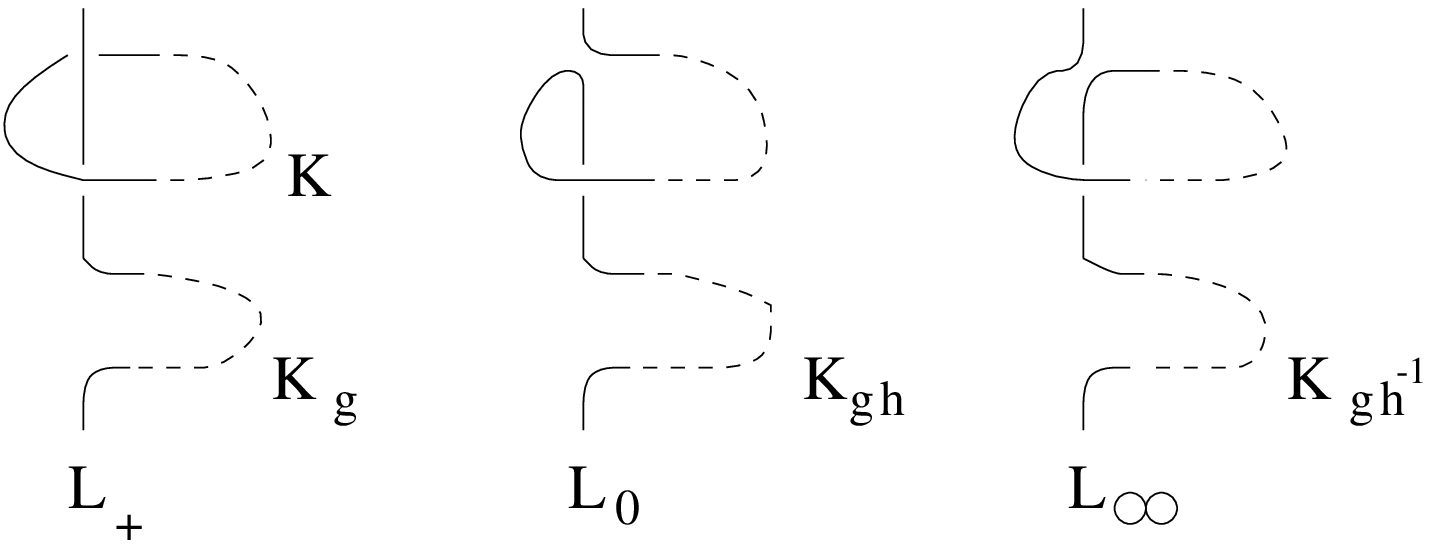}{3.5cm}}
\centerline{Figure 4}
\end{figure}

Therefore $\rho_0(L_+ +L_0+L_{\infty})=g\otimes \hat{\xi}([K][L''])+
gh\otimes \hat{\xi}([L''])+gh^{-1}\otimes \hat{\xi}([L''])+J=
(g\otimes \hat{\xi}([K])+ gh\otimes 1+gh^{-1}\otimes 1 +J)\cdot
(1\otimes \hat{\xi}([L''])+ J)=(g\otimes (-[h])+
gh\otimes 1+gh^{-1}\otimes 1+J)(1 \otimes \hat{\xi}([L''])+J)=0$
in $(RG\otimes_R {\cal S}(G;R))/J.$
\end{enumerate}
\end{proof}

Observe that if $L=K_g\cup L',$ where $L'$ is a link $F\times I,$ then
$\rho_0(L\cup \bigcirc + 2L)=\rho_0(L\cup \bigcirc)+ 2\rho_0(L)=
g\otimes \hat{\xi}([L'\cup \bigcirc])+2g\otimes \hat{\xi}([L'])+ J=
g\otimes \hat{\xi}(-2[L'])+2g\otimes \hat{\xi}([L'])+J= 0.$

The above observation and Lemma \ref{5.5} imply that
$\rho_0({\cal S}^{rel})=0,$ and hence, there is a homomorphism
$\rho_1:{\cal S}^{rel}_{2,\infty}(F\times I; R,-1)\to (RG\otimes_R
{\cal S}(G;R))/J,$
such that $\rho_1([L])=g\otimes \hat{\xi}([L'])+J$ for any relative framed 
link $L=K_g\cup L'.$

\begin{lemma}\label{5.6}\ \\
There is a homomorphism $\mu: (RG\otimes_R {\cal S}(G;R))/J \to H_R(G),$ such
that $\mu(g\otimes 1 + J)=[g]\in H_R(G),$ for any $g\in G.$
\end{lemma}

\begin{proof}
Let $\mu_0: RG\otimes {\cal S}(G;R)\to H_R(G)$ be a homomorphism of 
$R$-modules $\mu_0(x\otimes y)=[x]\cdot \psi(y),$ for any $x\in RG, y\in 
{\cal S}(G;R)$ \footnote {The homomorphism $\psi$ was introduced in Theorem 
\ref{3.2}.}.

Let $x_1,x_2\in RG, y_1,y_2\in {\cal S}(G;R).$ Then
$\mu_0((x_1\otimes y_1)(x_2\otimes y_2))=\mu_0((x_1x_2)\otimes (y_1y_2))=
[x_1x_2]\psi(y_1y_2)=[x_1][x_2]\psi(y_1)\psi(y_2).$ Since $\psi(y_1)\in TH_R(G)
\subset Center\ H_R(G),$ we have $\mu_0((x_1\otimes y_1)(x_2\otimes y_2))=
x_1\psi(y_1)x_2\psi(y_2)=\mu_0(x_1\otimes y_1)\mu_0(x_2\otimes y_2).$
Therefore $\mu_0$ is a homomorphism of $R$-algebras. Moreover,
$\mu_0(g\otimes [h]- gh\otimes 1 -gh^{-1}\otimes 1)=[g]([h]+[h^{-1}])-[gh]-
[gh^{-1}]=0$ for any $g,h\in G.$
Hence $\mu_0$ yields $\mu:(RG\otimes_R {\cal S}(G;R))/J \to H_R(G),$
such that $\mu(g\otimes 1+J)=[g].$
\end{proof}

Notice that $\mu\circ\rho_1\circ f_1:H_R(G)\to H_R(G).$
Let $g\in G.$ Then $\mu\circ\rho_1\circ f_1([g])=\mu\circ\rho_1([K_g])=
\mu(g\otimes 1+J)=[g].$ Since $\{[g]\}_{g\in G}$ generate $H_R(G),$
$\mu\circ\rho_1\circ f_1= id_{H_R(G)}.$ Therefore $f_1$ is a monomorphism,
and the proof of Theorem \ref{5.3} has been completed.
\end{proof}

Under some additional conditions we can strengthen the statement of Theorem 
\ref{5.3} in the following way.

Suppose that $A^{\pm 1}\in R$ and there is an involution $\tau_0$
on $R$ such that $\tau_0(A)=A^{-1}.$ For example
\begin{itemize}
\item $A=-1$ and $\tau_0: R\to R,\ \tau_0=id;$ \quad or
\item $R=R_0[A^{\pm 1}]$ for some ring $R_0$ and $\tau_0:R\to R,\ 
\tau_{0|R_0}=id_{R_0},\ \tau_0(A)=A^{-1}.$
\end{itemize}

Then we can define an involution $\tau$ on ${\cal S}^{rel}_{2,\infty}
(F\times I; R,A)$ in the following way.

Let $\tau_1: R{\cal L}_{fr}^{rel}(F\times I)\to R{\cal L}_{fr}^{rel}
(F\times I)$ be an additive function such that $\tau_1(rL)=\tau_0(r)h(L),$
where $h: F\times I\to F\times I$ is a homeomorphism $h(x,t)=(h,1-t),$ 
$x\in F, t\in [0,1].$ Then one can easily check that
$\tau_1({\cal S}_{fr}^{rel})={\cal S}_{fr}^{rel}.$ Therefore $\tau_1$
yields an additive function $\tau: {\cal S}^{rel}_{2,\infty}(F\times
I; R,A) \to {\cal S}^{rel}_{2,\infty}(F\times I; R,A).$ One can
further check that the following fact holds.

\begin{fact}\label{5.7}\ \\
$\tau$ is an involution on ${\cal S}^{rel}_{2,\infty}(F\times I;
R,A),$ i.e. $\tau$ is an anti-isomorphism of ${\cal S}^{rel}_
{2,\infty}(F\times I;R,A)$ and $\tau\circ\tau = id.$
Moreover, the following diagram commutes:\\

\begin{tabular}{cccc}
$R\hspace*{0.5cm}\times$ & ${\cal S}^{rel}_{2,\infty}(F\times I; R,A)$
& $\to$ & ${\cal S}^{rel}_{2,\infty}(F\times I; R,A)$\\
$\downarrow \tau_0$\hspace*{0.5cm} & $\downarrow \tau$& & $\downarrow \tau$\\
$R\hspace*{0.5cm}\times$ & ${\cal S}^{rel}_{2,\infty}(F\times I; R,A)$
& $\to$ & ${\cal S}^{rel}_{2,\infty}(F\times I; R,A)$
\end{tabular}
\end{fact}

Recall that we defined an involution $\imath$ on the algebra $H_R(G)$
at the beginning of Section 3.

\begin{theorem}\label{5.8}\ \\
Let $R$ be a commutative ring with the trivial involution $\tau_0=id_R.$
Then for any surface $F$ there exists a preserving involutions
isomorphism between $H_R(\pi_1(F))$ and 
${\cal S}^{rel}_{2,\infty}(F\times I; R,-1).$
\end{theorem}

\begin{proof}
Notice that for any $g\in \pi_1(F,p)$
$f_1(\imath([g]))=f_1([g^{-1}])=[K_{g^{-1}}]=\tau([K_g]).$ 
Therefore $f_1$ is an isomorphism of $R$-algebras which preserves involutions.
\end{proof}

We end this section with a theorem which can be considered as a generalization
of Theorem \ref{3.2}. Notice that if ${1\over 2}\in R$ then Theorem \ref{3.2}
(together with Theorems \ref{5.3} and \ref{2.8}) implies that there exist an 
injection ${\cal S}_{2,\infty}(F\times I;R,A)\to 
{\cal S}^{rel}_{2,\infty}(F\times I;R,A)$ for $A=-1.$ This fact turns
out to be true for any $A,$ and it has a simple topological
interpretation.

In order to introduce it consider a surface $F$ with a boundary and
a disc $D\subset F$ such that $D\cap \partial F$ is an arc in
$\partial F.$ Let $\gamma: [0,1]\to D\subset F$ be an arc and let
$K_0$ be the trivial special framed arc in $F\times I$ associated with
$\gamma,$ that is, $K_0: [0,1] \times I\to F\times I,\
K(s,t)=(\gamma(s),1-t).$ Since any framed link $L$ in $F\times I$ can
be pushed inside $(F\setminus D)\times I$ in a unique way (up to an
ambient isotopy in $(F\setminus D)\times I$), there is a well defined
operation of making a disjoint union of $L$ with $K_0.$ 

\begin{theorem}\label{5.9}\ \\
Let $F$ be a surface, $\partial F\ne \emptyset,$ and let
a disjoint union, $L\cup K_0\subset F\times I,$ be defined as above
for any framed link $L\subset F\times I.$ Then
there exists a homomorphism of algebras $\eta: {\cal S}_{2,\infty}
(F\times I; R,A)\to {\cal S}^{rel}_{2,\infty}(F\times I; R,A)$ such that
$\eta([L])=[L\cup K_0],$ for any framed link $L$ in $F\times I.$

\noindent Moreover, if $A^2+A^{-2}$ is invertible in $R$ then
$\eta: {\cal S}_{2,\infty}(F\times I; R,A)\to
{\cal S}^{rel}_{2,\infty}(F\times I; R,A)$ is a monomorphism.
\end{theorem}

\begin{proof}
Let $\eta_0: R{\cal L}_{fr}(F\times I)\to R{\cal L}^{rel}_{fr}(F\times
I)$ be
a homomorphism of $R$-modules such that $\eta_0(L)=L\cup K_0$ for any
link $L$ in $F\times I.$ Notice that this homomorphism 
induces a homomorphism of skein modules $\eta: {\cal S}_{2,\infty}
(F\times I;R,A)\to {\cal S}_{2,\infty}^{rel}(F\times I;R,A).$ Moreover,
$\eta$ is a homomorphism of rings.

Suppose now that $A^2+A^{-2}$ is invertible in $R.$ We will show that
$\eta$ is a monomorphism. 
Let $K_1$ be a framed arc in $D\times I\subset F\times I$ such that
$K_0\cap K_1=\gamma([0,1])\times \{0\}\cup \gamma([0,1])\times \{1\}$ and
$K_0\cup K_1,$ after pushing it inside the interior of $D\times I,$ is
a trivial framed knot in $F\times I.$
Then for any relative framed link $L$ in $F\times I,$
$L\cup K_1$ is a framed link which can be isotoped into the interior
of $F\times I.$ Moreover, we have a homomorphism
of modules $\nu:{\cal S}^{rel}_{2,\infty}(F\times I;R,A)\to {\cal
S}_{2,\infty}(F\times I; R,A)$ such that $\nu([L])=[L\cup K_1].$
Notice that $\nu\circ\eta :{\cal S}_{2,\infty}(F\times I;R,A)\to
{\cal S}_{2,\infty}(F\times I;R,A),\ \nu\circ\eta([L])=
\nu([L\cup K_0])=[L\cup K_0\cup K_1]=[L\cup \bigcirc]= -(A^2+A^{-2})[L].$
Therefore if $A^2+A^{-2}$ is invertible in $R$ then $\nu$ is a monomorphism.
\end{proof}

\section{Estimating minimal numbers of generators of skein algebras; Part 1}

In this section we will find minimal numbers of generators of skein algebras
of surfaces with coefficients in $Z[A^{\pm 1}].$ We will also calculate
minimal numbers of generators of skein algebras of abelian and non-abelian
free groups.

We start with a fact proven in [B-3] and concerning algebras
${\cal S}_{2,\infty}(F\times I; \Z[A^{\pm 1}],A)$ discussed in Section 4.

Let $F_{g,n}$ be a surface of genus $g$ with $n$ boundary components. If
$n\geq 1$ then we can present $F_{g,n}$ a disc with $2g+n-1$ handles
numbered as shown in Fig. 5.

\begin{figure}
 \centerline{\rysunek{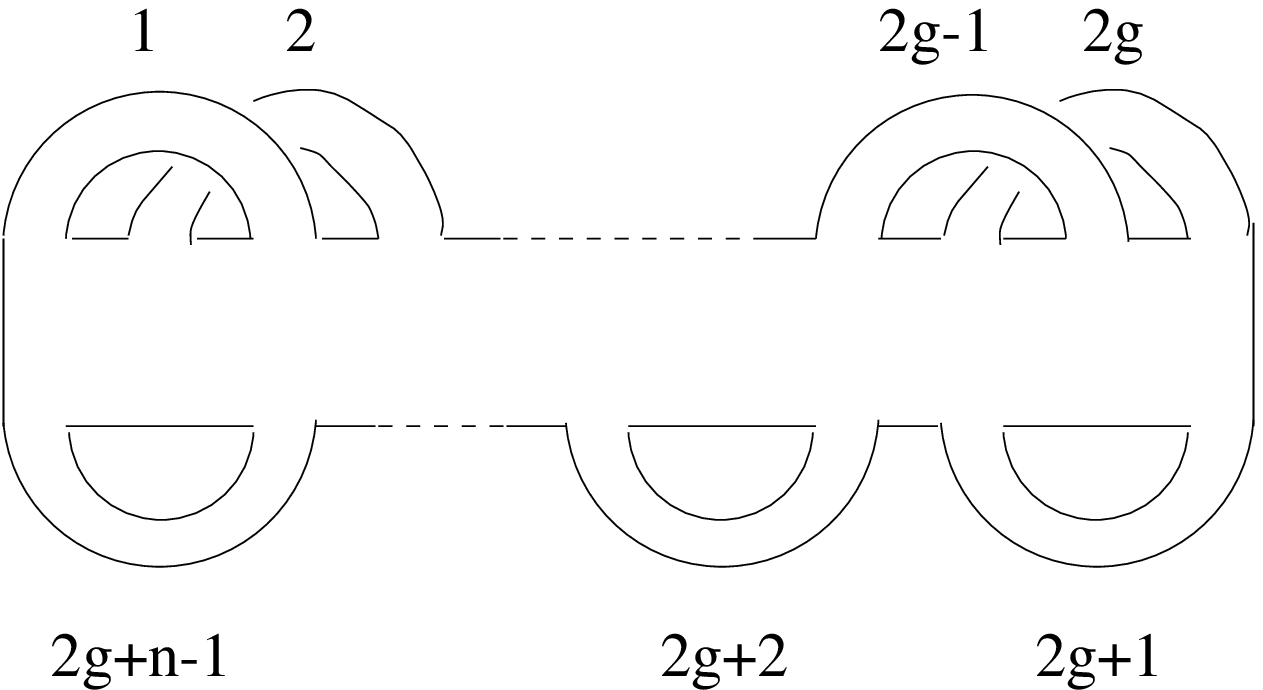}{3.2cm}}
 \centerline{Figure 5: Surface $F_{g,n}$}
\end{figure}

Let $S\subset \{1,2,...,2g+n-1\}.$ The standard knot of type $S$ in $F_{g,n},
n\geq 1,$ denoted by $K_S,$ is a knot which satisfies the following conditions:
\begin{enumerate}
\item if $i\in S$ then $K_S$ meets $i$-th handle exactly once.
\item if $i\not\in S$ then $K_S$ does not meet $i$-th handle.
\item if $K_S$ meets two overlapping handles (i.e. $2i+1,2i+2\in S,$ for some
$i\in \{0,1,...,g-1\}$) then it does it in a way shown in Fig. 6.
\end{enumerate}

\begin{figure}
 \centerline{\rysunek{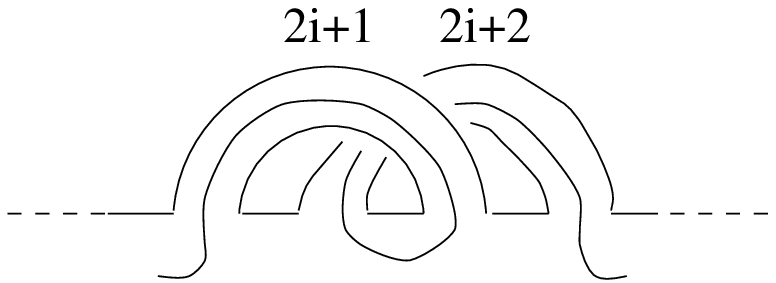}{2.3cm}}
 \centerline{Figure 6}
\end{figure}

We assume that $F_{g,1}\subset F_{g,0}$ and that the standard knots, $K_S,$ 
in $F_{g,1}$ are also the standard knots in $F_{g,0}.$ Therefore each
surface $F$ has $2^{rank\ H_1(F)}$ standard knots.
 
\begin{theorem}[\cite{B-3}]\label{6.1}\ \\ 
${\cal S}_{2,\infty}(F_{g,n}\times I;\Z[A^{\pm 1}],A)$ is a 
$\Z[A^{\pm 1}]$-algebra generated by $2^N-1$ elements $[K_S],\ N=rank
H_1(F),$ for all possible subsets $S\subset \{1,2,...,N\},\ S\ne\emptyset.$
\end{theorem}

Using Bullock's proof of Theorem \ref{6.1}, we can prove the following
theorem.

\begin{theorem}\label{6.2}\ \\
Let $M$ be a compact manifold of a Heegaard genus $g.$
Then there is a link $L$ in $M$ of $N=2^g-1$ components, $K_1, K_2,...,
K_N,$ such that for any family of disjoint regular neighborhoods 
$K_i\subset V_i\simeq S^1\times D^2,\ i=1,...,N,$ the embedding 
$i: V_1\cup V_2\cup ... \cup V_N \to M$ induces
an epimorphism of $R$-modules $i_*: {\cal S}_{2,\infty}(V_1;R,A)\otimes
{\cal S}_{2,\infty}(V_2;R,A)\otimes ... \otimes {\cal S}_{2,\infty}(V_N;R,A)
\to {\cal S}_{2,\infty}(M;R,A).$
\end{theorem}

\begin{proof}Let us consider any linear order, $<,$ on the set of non-empty 
subsets of $\{1,2,...,N\},$ where $N=rank H_1(F).$ Bullock's proof of 
Theorem \ref{6.1} can be modified to give a stronger version of 
Theorem \ref{6.1}. Namely, 
${\cal S}_{2,\infty}(F_{g,n}\times I;\Z[A^{\pm 1}],A)$
is a $\Z[A^{\pm 1}]$-module generated by monomials $[K_{S_1}]\cdot [K_{S_2}]
\cdot ... \cdot [K_{S_k}]$ such that $S_1\leq S_2\leq ...\leq S_k.$

Suppose that $M$ has a Heegaard splitting $M=H_1\cup H_2$ of genus $g.$
Then $H_1$ is homeomorphic to $F_{0,g+1}\times [0,1].$ 
Let $L$ be a link in $H_1\simeq F_{0,g+1}\times [0,1]$ composed of 
$2^g-1$ knots $K_S,$ where $S\subset \{1,...,g\},S\ne \emptyset,$ in such 
a way that $K_{S_1}$ lies `above' $K_{S_2}$ in $H_1 \simeq 
F_{0,g+1}\times I,$ iff $S_1>S_2.$
Notice that then the embedding $j: V_1\cup V_2\cup ... \cup V_N \to H_1$ 
induces an epimorphism of $R$-modules 
$$j_*: {\cal S}_{2,\infty}(V_1;R,A)\otimes {\cal S}_{2,\infty}(V_2;R,A)\otimes
 ... \otimes {\cal S}_{2,\infty}(V_N;R,A)\to {\cal S}_{2,\infty}
(H_1;R,A).$$ Moreover, since every link in $M$ can be pushed inside 
$H_1$ we have an epimorphism 
${\cal S}_{2,\infty}(H_1;\Z[A^{\pm 1}],A)\to
{\cal S}_{2,\infty}(M\times I;\Z[A^{\pm 1}],A).$ The composition of this
epimorphism with $j_*$ gives the required epimorphism $i_*.$
\end{proof}

Theorem \ref{6.1} implies the following theorem.

\begin{theorem}\label{6.3}\ \\
\begin{enumerate}
\item Let $F_n$ denote a free group, $F_n=<g_1,g_2,...,g_n>.$
Then for any ring $R$ the skein algebra ${\cal S}(F_n; R)$ is generated by 
$2^n-1$ elements of the form $[g_{i_1}g_{i_2}...g_{i_k}],$ where 
$1\leq i_1<i_2<...<i_k\leq n$ and $k\geq 1.$ 
\item For any ring $R$ the skein algebra ${\cal S}(\Z^n; R)$ is generated by 
$2^n-1$ elements $[(\epsilon_1,...,\epsilon_n)],$ where 
$(\epsilon_1,...,\epsilon_n)\in \Z^n,$
$\epsilon_1,\epsilon_2,...,\epsilon_n\in \{0,1\}$ and $\epsilon_i\ne 0$ for 
at least one $i.$
\item If $R=\Z$ then the numbers of generators given above are minimal, i.e.
neither ${\cal S}(F_n; \Z)$ nor ${\cal S}(\Z^n; \Z)$ can have fewer generators
than $2^n-1.$ \footnote{In fact we prove the following more general fact: 
If $R$ is a ring such that $\Z_2$ can be an epimorphic image of $R$ then 
neither ${\cal S}(F_n; R)$ nor ${\cal S}(\Z^n; R)$ can have fewer generators 
than $2^n-1.$}
\end{enumerate}
\end{theorem}

\begin{proof}
\begin{enumerate}

\item
We can consider $R$ as a $\Z[A^{\pm 1}]$-algebra determined by the 
homomorphism $f:\Z[A^{\pm 1}]\to R, f(1)=1, f(A)=-1.$ Theorem \ref{6.1} and
Fact \ref{2.7} implies that the algebra ${\cal S}_{2,\infty}(F_{0,n+1}\times I
;R,-1)$
is generated by $[K_S],\ S\subset \{1,2,...,n\}, S\ne \emptyset.$
Let us fix an isomorphism between $F_n=<g_1,g_2,...,g_n>$ and 
$\pi_1(F_{0,n+1})$ in such a way that $g_i$ corresponds to the homotopy
class of $K_{\{i\}}.$
Then, by Theorem \ref{2.8}, ${\cal S}_{2,\infty}(F_{0,n+1}\times I;R,-1)$ 
is isomorphic to ${\cal S}(F_n; R)$ and the generators $[K_S],$ for 
$S=\{i_1,i_2,...,i_k\}\subset \{1,2,...,n\},$ of the skein algebra
${\cal S}_{2,\infty}(F_{0,n+1}\times I;R,-1)$ correspond to elements 
$(-1)^{k}[g_{i_1}g_{i_2}...g_{i_k}]\in {\cal S}(F_n;R).$ 

\item
Let $f:F_n=<g_1,...,g_n> \to \Z^n$ be an epimorphism given by $f(g_i)=
(0,0,...,1,...0),$ where the single $1$ stands at the $i$-th place. By Fact 
\ref{2.6}(3), $f$ yields an epimorphism $f_*:{\cal S}(F_n; R)\to {\cal S}
(\Z^n; R).$
Notice that $f_*([g_{i_1}g_{i_2}\ldots g_{i_k}])=
[(\epsilon_1,...,\epsilon_n)],$ where $\epsilon_1,\epsilon_2,...,\epsilon_n
\in \{0,1\}$ and $\epsilon_j=1,$ for $j=1,2,...,n,$ iff $j$ is equal to
one of the numbers $i_1,i_2,...,i_k.$ 
Therefore ${\cal S}(\Z^n; R)$ is generated by all elements of the form 
$[(\epsilon_1,...,\epsilon_n)],$ where $(\epsilon_1,...,\epsilon_n)\in \Z^n,$
$\epsilon_i\in \{0,1\}$ and at least one $\epsilon_i$ is equal
to $1.$

\item
Let $h:\Z^n\to \Z_2^n$ be a natural projection sending $(x_1,...,x_n)\in \Z^n$
to \ \\
$(x_1 \mod 2,...,x_n \mod 2)\in \Z_2^n.$ Let $f:F_n\to \Z^n$ be defined as
in the previous paragraph.
Since $h_*:{\cal S}(\Z^n; \Z)\to {\cal S}(\Z_2^n; \Z)$ and 
$h_*\circ f_*:{\cal S}(F_n; \Z)\to {\cal S}(\Z_2^n; \Z)$ are epimorphisms, it
is enough to show that ${\cal S}(\Z_2^n; \Z)$ requires at least $2^n-1$
generators. Suppose that ${\cal S}(\Z_2^n; \Z)$ can be generated by k elements.
By Fact \ref{2.6}(4), ${\cal S}(\Z_2^n; \Z_2)=
{\cal S}(\Z_2^n; \Z)\otimes_\Z \Z_2.$
Therefore ${\cal S}(\Z_2^n; \Z_2)$ can also be generated by $k$ elements.
Let $g,h\in \Z_2^n.$ Then $gh=gh^{-1}$ and $[g]\otimes [h]=[gh]+[gh^{-1}]=
2[gh]=0$ in ${\cal S}(\Z_2^n; \Z_2)$ (We use the multiplicative notation for
the addition in $\Z_2^n$). Moreover, $[e]=0$ in ${\cal S}(\Z_2^n; \Z_2)$
($e$ is the identity in $\Z_2^n$). Therefore ${\cal S}(\Z_2^n; \Z_2)$ is 
isomorphic to the ring $R=\Z_2[{x_g};\ {g\in \Z_2^n}]/(x_e,x_gx_h),$ where 
$\Z_2[{x_g};\ {g\in \Z_2^n}]$ means the ring of polynomials in $2^n$ variables
$x_g,\ g\in \Z_2^n.$
We will show that the ring $R$ requires at least $2^n-1$ generators.
Notice that for any $x,y\in R$ $xy=0$ or $xy=1.$ Therefore $R$ is a
$\Z_2$-linear space of dimension $2^n$ with a basis 
$\{x_g\}_{g\in \Z^n_2\setminus \{e\}}\cup \{1\}.$
If $y_1,y_2,...,y_k$ generate $R$ than $1,y_1,y_2,...,y_k$ span $R$ as a 
linear space (because $y_iy_j=0$ or $y_iy_j=1$). Therefore $k\geq 2^n-1.$
\end{enumerate}
\end{proof}

Notice that Fact \ref{2.6}(3) implies the following corollary to Theorem
\ref{6.3}.

\begin{corollary}\label{6.4}\ \\
Let $G$ be a group generated by $g_1,g_2,...,g_n\in G.$ Then for any
ring $R$ the skein algebra ${\cal S}(G;R)$ is generated by elements of the form
$[g_{i_1}g_{i_2}...g_{i_k}],$ where $1\leq i_1<i_2<...<i_k\leq n$ and
$k\geq 1.$ 
\end{corollary}

Using Theorem \ref{6.3} we can prove the following result.

\begin{proposition}\label{6.5}\ \\
The algebra ${\cal S}_{2,\infty}(F_{g,n}\times I;\Z[A^{\pm 1}],A)$ cannot 
be generated by fewer than $2^N-1$ elements, where $N=rank H_1(F).$
\end{proposition}

\begin{proof}
By Fact \ref{2.7} $${\cal S}_{2,\infty}(F_{g,n}\times I;\Z,-1)\simeq
{\cal S}_{2,\infty}(F_{g,n}\times I;\Z[A^{\pm 1}],A)\otimes_{\Z[A^{\pm 1}]} 
\Z,$$ where $\Z$ is considered as a $\Z[A^{\pm 1}]$-algebra via homomorphism
$h:\Z[A^{\pm 1}]\to \Z,\ h(A)=-1.$ Let $k$ be the minimal number of generators
of ${\cal S}_{2,\infty}(F_{g,n}\times I;\Z[A^{\pm 1}],A).$ Then 
${\cal S}_{2,\infty}(F_{g,n}\times I;\Z, -1)\simeq {\cal S}
(\pi_1(F_{g,n}); \Z)$ can also be generated by $k$ elements. 
Let $f:\pi_1(F_{g,n})\to \Z^N$ be the abelianizing homomorphism. Then
$f_*: {\cal S}(\pi_1(F_{g,n}); \Z)\to {\cal S}(\Z^N; \Z)$ is onto and
by Theorem \ref{6.3}(3) $k\geq 2^N-1.$
\end{proof}

\section{Characters of $Sl_2(\C)$-representations of groups}

In this section we combine the results of previous sections with the results
of [B-H] and apply them to a study of connections between skein
algebras and $SL_2(\C)$-character varieties.
We start with a brief presentation of the approach to 
$SL_2(\C)$-representations as discussed in [B-H]. 
For a more general theory see [L-M].

For any finitely generated group $G$ one can assign a commutative 
{\it universal
representation $\C$-algebra,} $A[G],$ and a {\it universal representation}
$$\rho_G: G\to SL_2(A[G]),$$ with the following universal property:

\noindent For any $\C$-algebra $A$ and any representation $\rho: G\to SL_2(A)$
there is a unique homomorphism of algebras $h_{\rho}: A[G]\to A$ which
induces a homomorphism of groups
$$SL_2(h_{\rho}):SL_2(A[G])\to SL_2(A)$$
such that the following diagram commutes:

\begin{center}
\begin{tabular}{ccc}
$G$ & $\stackrel{\rho_G}{\longrightarrow}$ & $SL_2(A[G])$\\
& $\stackrel{\rho}{\searrow}$ & $\Big\downarrow
\scriptstyle{SL_2(h_{\rho})}$\\
& & \hspace*{-0.5in} $SL_2(A)$\\
\end{tabular}
\end{center}

\noindent This universal property uniquely determines $A[G]$ up to
an isomorphism.
One can think of $A[G]$ as an (unreduced) coordinate ring of a scheme
whose points, that is homomorphisms $A[G]\to \C,$ are in a natural
bijection with $SL_2(\C)$-representations of $G.$ $GL_2(\C)$ acts
on these representations by conjugation. A standard argument from
invariant theory states that there is a categorical quotient of
this action, $A[G]^{GL_2(\C)},$ called {\it a universal character
ring of $G$.} Its $\C$-points represent semisimple 
$SL_2(\C)$-representations of $G,$ or equivalently
$SL_2(\C)$-characters of $G.$ $A[G]^{GL_2(\C)}$ is the most
important object in the study of $SL_2$-representations of $G$ since it
encodes all information necessary for their understanding
(or classification).
One of the most important results of [B-H], Proposition 9.1 (as well as
$A^*$.9.11), gives a simple algebraic description of $A[G]^{GL_2(\C)}$
by showing that $A[G]^{GL_2(\C)}$ is isomorphic to the algebra $TH_{\C}(G)$
introduced in Section 3. The isomorphism carries
$[g]+[g^{-1}]\in TH_{\C}(G)\subset
H_{\C}(G),$ for $g\in G,$ to $Tr(\rho_G(g)).$ This result and 
Theorem \ref{3.2} imply the following theorem.

\begin{theorem}\label{7.1}\ \\
There is an isomorphism of $\C$-algebras $\phi: {\cal S}(G;R)\to
A[G]^{GL_2(\C)}$ such that $\phi([g])=Tr(\rho_G(g)),$ for any $g\in
G.$
\end{theorem}

The above theorem implies in particular the following corollary.

\begin{corollary}\label{7.2}\ \\
\begin{enumerate}
\item Let $G$ be any finitely generated group and $\chi:G\to \C$ be an
$Sl_2(\C)$-character of $G,$ that is, the trace of an
$SL_2(\C)$-representation of $G.$
Then there is a unique homomorphism of $\C$-algebras
$h_{\chi}:{\cal S}(G; \C)\to \C$ such that $h_{\chi}([g])=\chi(g).$
If $\chi$ and $\chi'$ are two different $Sl_2(\C)$-characters then
$h_{\chi}\ne h_{\chi'}.$
\item Conversely, every homomorphism
$h:{\cal S}(G;R)\to \C$ is induced by an $SL_2(\C)$-character of $G,$ i.e.
$h=h_{\chi},$ for some $\chi:G\to \C.$
\end{enumerate}
\end{corollary}

{\bf Remarks}\\
\begin{enumerate}
\item Doug Bullock was the first to observe that the well-known equality
(see e.g [Vo], [F-K], [Ho])
$$\chi(a)\chi(b)=\chi(ab)+\chi(ab^{-1}),$$ satisfied for
any $SL_2(\C)$-character $\chi$ and $g,h\in G$ is very similar to
the skein relation
$$[L_+]=A[L_0]+A^{-1}[L_{\infty}] {\rm\ in\ } {\cal S}_{2,\infty}(M;R,A)
{\rm\ for\ } A=-1.$$ As an immediate consequence of this observation he 
got Corollary \ref{7.2}(1); compare [B-1], [B-2], [P-S-1]. This is the
easier part of Corollary \ref{7.2}. He also proved independently from
us Corollary \ref{7.2}(2), [B-4].
\item The statement of Theorem \ref{7.1} is much stronger than
Corollary \ref{7.2}. For example $A[G]^{GL_2(\C)}$ may have nilpotent
elements, which obviously will be undetected by homomorphisms
$A[G]^{GL_2(\C)}\to \C.$ For more information on $A[G]^{GL_2(\C)},$
and a discussion of possible nilpotent elements in this ring see
[B-H]. We will come back to the question of nilpotents later in this
section.
\item The second author proved (see [Si]) that Corollary \ref{7.2}
is  true not only for the field of complex numbers but for all algebraically 
closed fields of characteristic $\ne 2.$ Moreover, it can be generalized to
fields which are not closed and to Dedekind domains.
\item The connections between skein modules and character varieties can be
generalized to other algebraic groups. In particular, an analogous
description of the universal $SL_n(\C)$-character ring of $\pi_1(M)$
in terms of links in $M$ is given in [S-F].
\end{enumerate}

Corollary \ref{7.2} can be formulated in a different way.
Let $X(G)$ denote the set of all $SL_2(\C)$-characters of $G.$
For any $g\in G$ there is a function $\tau_g:X(G)\to \C$ given by
$\tau_g(\chi)=\chi(g).$ It is known that $X(G)$ has a natural
structure of an algebraic set determined by the requirement that each
$\tau_g$ is a regular function on $X(G);$ see for example [C-S], [Bau].
The algebraic set $X(G)$ is called {\it the $SL_2(\C)$-character
variety of $G.$} Its coordinate ring, $\C[X(G)],$ is isomorphic to
$(A[G]^{GL_2(\C)})/\sqrt 0,$ and therefore to ${\cal S}(G;\C)/\sqrt
0,$ where $\sqrt 0$ denotes the radical of the ring, that is the ideal
composed of all nilpotent elements. The isomorphism between
${\cal S}(G; \C)/\sqrt 0$ and $\C[X(G)]$ carries $[g]$ to $\tau_g.$
In this context it is natural to ask when ${\cal S}(G;\C)$ has no
nilpotent elements and therefore ${\cal S}(G; \C)\simeq \C[X(G)].$ The next
theorem gives a list of all classes of groups for which we know a positive
answer to this question.

\begin{theorem}\label{7.3}\ \\
If one of the following conditions is satisfied
\begin{enumerate}
\item $G$ is finite
\item $G$ is abelian
\item $G$ is free non-abelian
\item $G$ is the fundamental group of a surface satisfying condition (1) or
(2) of Theorem \ref{4.7}
\item $G$ is the knot group of a $2$-bridge knot
\end{enumerate}
then ${\cal S}(G;\C)$ does not have nilpotent elements and, therefore,
${\cal S}(G;\C)\simeq \C[X(G)].$
\end{theorem}

\begin{proof}\ \\
(1) If $G$ is finite then by Maschke's Theorem $\C G$ is a semisimple algebra.
By Wedderburn-Artin theory (see e.g. [C-R],[Hu]) any homomorphic image of a 
semisimple, artinian algebra is also semisimple and artinian. Therefore 
$\C G/I,$ where $I$ is generated by $g(h+h^{-1})-(h+h^{-1})g,\ g,h\in G,$ is
semisimple and artinian. Hence, by Wedderburn-Artin Theorem there exists an
isomorphism $\phi:\C G/I\to M_{n_1}(\C)\times\ldots\times M_{n_k}(\C),$
where $M_{n_i}(\C)$ is the algebra of $n_i\times n_i$ matrices with complex 
coefficients. Since $TH_\C(G)$ is generated by the elements
$g+g^{-1}+I\in \C G/I$ and these elements are in the center of $\C G/I,$
$\phi(TH_\C(G))\subset Center(M_{n_1}(\C)\times\ldots\times M_{n_k}(\C))=
\C^k.$
Therefore $\phi(TH_\C(G))$ has no nilpotent elements and, finally, $TH_\C(G)$
is nilpotent free, because $\phi$ is an isomorphism.

(2) It has been proved in [P-S-1] (Theorem 3.3).\ \\

(4) It follows immediately from Theorem \ref{4.7}\ \\

(3) This is a special case of (4).\ \\

(5) The Burde's and Zieschang's book [B-Z] is a good source of information
about $2$-bridge knots. Our proof will base on results of L. Thang [Th] on
$Sl_2(\C)$-representations of $2$-bridge knots
\footnote{These representations were also studied in [Ri] and [B-H].}.
Let $K$ be a $2$-bridge knot. Then its group $G=\pi_1(S^3\setminus K)$ has a 
presentation $G=<a,b | wa=bw>,$ where $w=a^{\epsilon_1}b^{\epsilon_n}
a^{\epsilon_2}b^{\epsilon_{n-1}}...a^{\epsilon_n}b^{\epsilon_1},
\epsilon_i=\pm 1.$ Let $F_2=<a,b>.$ By Fact \ref{2.6}(3) the natural 
epimorphism $\pi:F_2\to G$ yields an epimorphism $\pi_*:
{\cal S}(F_2;\C)\to {\cal S}(G;\C).$ Let $t_1=[a],\ t_2=[ab],\ t_3=[b],\ 
t_1,t_2,t_3\in {\cal S}(F_2;\C).$ Then by Corollary \ref{4.5}
${\cal S}(F_2; \C)=\C[t_1,t_2,t_3].$ Since $[w], [bwa^{-1}]\in 
{\cal S}(F_2,\C),\ [w]=P_w(t_1,t_2,t_3), [bwa^{-1}]=P_{bwa^{-1}}
(t_1,t_2,t_3),$ where $P_w, P_{bwa^{-1}}$ are polynomials in variables 
$t_1,t_2,t_3.$ Then 
\begin{equation}
P_w-P_{bwa^{-1}}\in Ker \pi_* \triangleleft \C[t_1,t_2,t_3] \label{e1}
\end{equation}
Moreover, $a$ and $b$ are conjugate in $G.$
Therefore $[a]=[b]$ in ${\cal S}(G; \C)$ and 
\begin{equation}
t_1-t_3\in Ker \pi_*. \label{e2}
\end{equation}
Since $[a]$ and $[ab]$ generate ${\cal S}(G;\C),$ $a$ and $ab$ distinguish 
all $Sl_2(\C)$-characters of $G.$ Hence, by identifying any $\chi\in X(G)$ 
with a pair $(t_1,t_2)\in \C^2,$ where $t_1=\chi(a),t_2=\chi(ab),$ we can 
consider $X(G)$ as a subset of $\C^2.$ Le Thang has shown in [Th] (Theorem 
3.3.1) that $X(G)$ is the zero set of a polynomial equation 
$P_w(t_1,t_2,t_1)-P_{bwa^{-1}}(t_1,t_2,t_1)=0.$ Moreover, he has shown (see 
[Th] Proposition 3.4.1) that
$P_w(t_1,t_2,t_1)-P_{bwa^{-1}}(t_1,t_2,t_1)=(t_1^2-t_2-2)\Phi(t_1,t_2),$
where $\Phi\in \C[t_1,t_2]$ has no multiple divisors. He also proved 
([Th] Lemma 3.3.6) that $\Phi(2,2)\ne 0.$ Since $t_1^2-t_2-2=0$ for 
$t_1=t_2=2,$ $t_1^2-t_2-2$ does not divide $\Phi(t_1,t_2).$ Moreover,
$t_1^2-t_2-2$ is an irreducible polynomial. Therefore 
$P_w(t_1,t_2,t_1)-P_{bwa^{-1}}(t_1,t_2,t_1)$ has no quadratic divisors.
This implies that $(P_w(t_1,t_2,t_1)-P_{bwa^{-1}}(t_1,t_2,t_1))\triangleleft
\C[t_1,t_2]$ is a radical ideal and 
\begin{equation}
\C[X(G)]=\C[t_1,t_2]/(P_w(t_1,t_2,t_1)-P_{bwa^{-1}}(t_1,t_2,t_1)) \label{e3}
\end{equation}
From the discussion preceding Theorem \ref{7.3} follows that there is
a homomorphism
\[\phi:{\cal S}(G; \C)\to \C[X(G)]=\C[t_1,t_2]/(P_w(t_1,t_2,t_1)-
P_{bwa^{-1}}(t_1,t_2,t_1))\]
such that $\phi([a])=t_1,\phi([ab])=t_2.$
Therefore 
\[\phi\circ\pi_*:\C[t_1,t_2,t_3]\to {\cal S}(G; \C)\to \C[X(G)]=
\C[t_1,t_2]/(P_w(t_1,t_2,t_1)-P_{bwa^{-1}}(t_1,t_2,t_1))\]
maps $t_1,t_3$ on $t_1$ and $t_2$ on $t_2.$ Hence $Ker \pi_*\subset
(t_1-t_3, P_w(t_1,t_2,t_1)-P_{bwa^{-1}}(t_1,t_2,t_1))=
(t_1-t_3, P_w(t_1,t_2,t_3)-P_{bwa^{-1}}(t_1,t_2,t_3))\triangleleft 
\C[t_1,t_2,t_3].$ But we have observed in (\ref{e1}) and (\ref{e2}) that 
$t_1-t_3\in Ker \pi_*$ and $P_w(t_1,t_2,t_3)-P_{bwa^{-1}}(t_1,t_2,t_3)\in 
Ker \pi_*.$
Therefore $Ker \pi_*=(P_w(t_1,t_2,t_1)-P_{bwa^{-1}}(t_1,t_2,t_1)).$
Hence ${\cal S}(G; \C)\simeq \C[t_1,t_2,t_3]/Ker \pi_*\simeq \C[X(G)]$ has 
no nilpotent elements.
\end{proof}

Of course Theorem \ref{4.7} implies much more than the statement of
Theorem \ref{7.3}(4).

\begin{corollary}\label{7.4}\ \\
Suppose that $F$ is either an orientable surface or
$F$ is an unorientable surface of an even, negative Euler characteristic.
Then the $Sl_2(\C)$-character variety $X(\pi_1(F))$ is an irreducible
affine algebraic set.
\end{corollary}

$Sl_2(\C)$-character varieties of surfaces were investigated earlier by
several authors. In particular, W. Goldman constructed a Poisson bracket
on coordinate rings of $Sl_2(\C)$-character varieties of closed surfaces
(See [Go]). D. Bullock, C. Frohman and J. Kania-Bartoszy\' nska showed
that the skein algebra ${\cal S}_{2,\infty}(F\times I;\C[A^{\pm 1}],A)$ can be 
considered as a quantization of $\C[X(\pi_1(F))]$ according to the
Goldman-Poisson bracket, [B-F-K].

Using our Theorem \ref{7.1} one can reformulate some of the results of 
Culler and Shalen stated in [C-S]. In particular, they proved that
if $M$ is a compact, oriented 3-manifold and an irreducible
component of the $Sl_2(\C)$-character variety of $\pi_1(M)$
has dimension greater or equal to $1$ then $M$ has an
incompressible surface not parallel to the boundary of $M.$
Bullock observed that the result of Culler and Shalen combined with the 
statement of Theorem \ref{7.1} implies the following corollary.

\begin{corollary}\label{7.5}\ \\
If the skein module, ${\cal S}_{2,\infty}(M;\C,-1),$ of a compact, 
oriented 3-manifold $M$ is infinite dimensional (as a $\C$-linear
space) then $M$ has an incompressible surface which is not boundary-parallel.
\end{corollary}

The above corollary gives a topological condition on links in $M$
(modulo skein relations) which implies the existence an incompressible
surface in $M.$ The only existing proof of Corollary \ref{7.5} is
very complicated. It would be very desirable to find a new, simpler and
purely topological proof of this corollary.

Corollary \ref{7.5} shows that ${\cal S}_{2,\infty}(M;\Z[A^{\pm 1}],A)$ 
carries a great amount of important information about the manifold $M.$
\footnote{ Notice that ${\cal S}_{2,\infty}(M;\C,-1)=
{\cal S}_{2,\infty}(M;\Z[A^{\pm 1}],A)\otimes_{\Z[A^{\pm 1}]} \C.$}

\section{Estimating minimal numbers of generators of skein algebras, Part 2}

In this section we are going to prove theorems similar to theorems of
Section 6, but concerning the case in which $A^2+A^{-2}$ is invertible in 
rings of coefficients of skein algebras.
We use the notation introduced in Section 6.

\begin{theorem}\label{8.1}\ \\
Let $A$ be an invertible element in a ring $R$ such that
$(A^2+A^{-2})^{-1}\in R.$ Then ${\cal S}_{2,\infty}(F_{g,n}\times I; R,A)$ 
is generated by $N+{N \choose 2}+{N \choose 3}$ elements $[K_S],$
where $N=rank H_1(F_{g,n})$ and $S\subset \{1,...,N\}, S\ne \emptyset, 
\card{S}\leq 3$ \footnote{$\card{S}$ denotes the number of elements of $S.$}.
\end{theorem}

\begin{proof}
Since any link in $F_{g,0}$ can be pushed into $F_{g,1}\subset F_{g,0},$
the embedding $i:F_{g,1}\to F_{g,0}$ induces an epimorphism
$i_*:{\cal S}_{2,\infty}(F_{g,1}\times I; R,A)\to {\cal S}_{2,\infty}
(F_{g,0}\times I; R,A).$ Hence, it is enough to prove Theorem \ref{8.1}
only for surfaces $F_{g,n},\ n\geq 1.$
 
Let us consider a homomorphism $r:\Z[A^{\pm 1}]\to R,\ r(A)=A\in R.$
By Fact \ref{2.7} there is an isomorphism
\[{\overline r}:{\cal S}_{2,\infty}(F_{g,n}\times I;R,A) \to
{\cal S}_{2,\infty}(F_{g,n}\times I;\Z[A^{\pm 1}],A)\otimes_{\Z[A^{\pm 1}]} 
R.\]
Therefore Theorem \ref{6.1} implies that 
${\cal S}_{2,\infty}(F_{g,n}\times I; R,A)$ is generated by the elements 
$[K_S], S\subset \{1,...,N\}.$ Hence it is enough to prove that the 
subalgebra $P\subset {\cal S}_{2,\infty}(F_{g,n}\times I; R,A)$ generated by 
$[K_S], \card{S}\leq 3,$ contains all elements $[K_S], S\subset 
\{1,2,...,N\},\ S\ne \emptyset.$

Suppose that $S_0=\{i_1,i_2,...,i_k\}\subset \{1,2,...,N\}, i_1<i_2<...<i_k,
k\geq 4,$ is a set of the smallest cardinality such that $[K_{S_0}]\not\in P.$
Let us consider knots $K_{\{i_1,i_3\}}$ and $K_{S_0\setminus \{i_1,i_3\}}.$ 
Then, depending on $S_0,$ the minimal number of intersections between 
$K_{\{i_1,i_3\}}$ and $K_{S_0\setminus \{i_1,i_3\}}$ is equal to $2,3$ or $4$.
Let us consider each of these possibilities.
\begin{enumerate}
\item The minimal number of intersections between $K_{\{i_1,i_3\}}$ and 
$K_{S_0\setminus \{i_1,i_3\}}$ is four.\\
One can check by straightforward, but long computations that

$[K_{\{i_1,i_3\}}][K_{S_0\setminus \{i_1,i_3\}}]=
(A^2+A^{-2})[K_{S_0}]+ 
A^4[K_{\{i_1,i_2\}}][K_{S_0\setminus\{i_1,i_2\}}]+$\\
$A^{-4}[K_{\{i_2,i_3\}}][K_{S_0\setminus \{i_2,i_3\}}]+
A^2[K_{\{i_1,i_2\}}][K_{\{i_3\}}][K_{S_0\setminus\{i_1,i_2,i_3\}}]+$\\
$A^2[K_{\{i_1\}}][K_{\{i_2\}}][K_{S_0\setminus \{i_1,i_2\}}]+
A^{-2}[K_{\{i_2\}}][K_{\{i_3\}}][K_{S_0\setminus \{i_2,i_3\}}]+$\\
$A^{-2}[K_{\{i_2,i_3\}}][K_{\{i_1\}}][K_{S_0\setminus \{i_1,i_2,i_3\}}]+
[K_{\{i_1\}}][K_{\{i_2\}}][K_{\{i_3\}}][K_{S_0\setminus\{i_1,i_2,i_3\}}]+$\\
$[K_{\{i_2\}}][K_{S_0\setminus \{i_2\}}]+
[K_{\{i_3\}}][K_{S_0\setminus \{i_3\}}]+$\\
$[K_{\{i_1,i_2,i_3\}}][K_{S_0\setminus \{i_1,i_2,i_3\}}]+
[K_{\{i_1\}}][K_{S_0\setminus \{i_1\}}].$\\

The above calculations for the knot $K_{S_0}$ placed in $F_{0,5},$
$S_0=\{1,2,3,4\},$ are shown in Fig. 7.

\begin{figure}
 \centerline{\rysunek{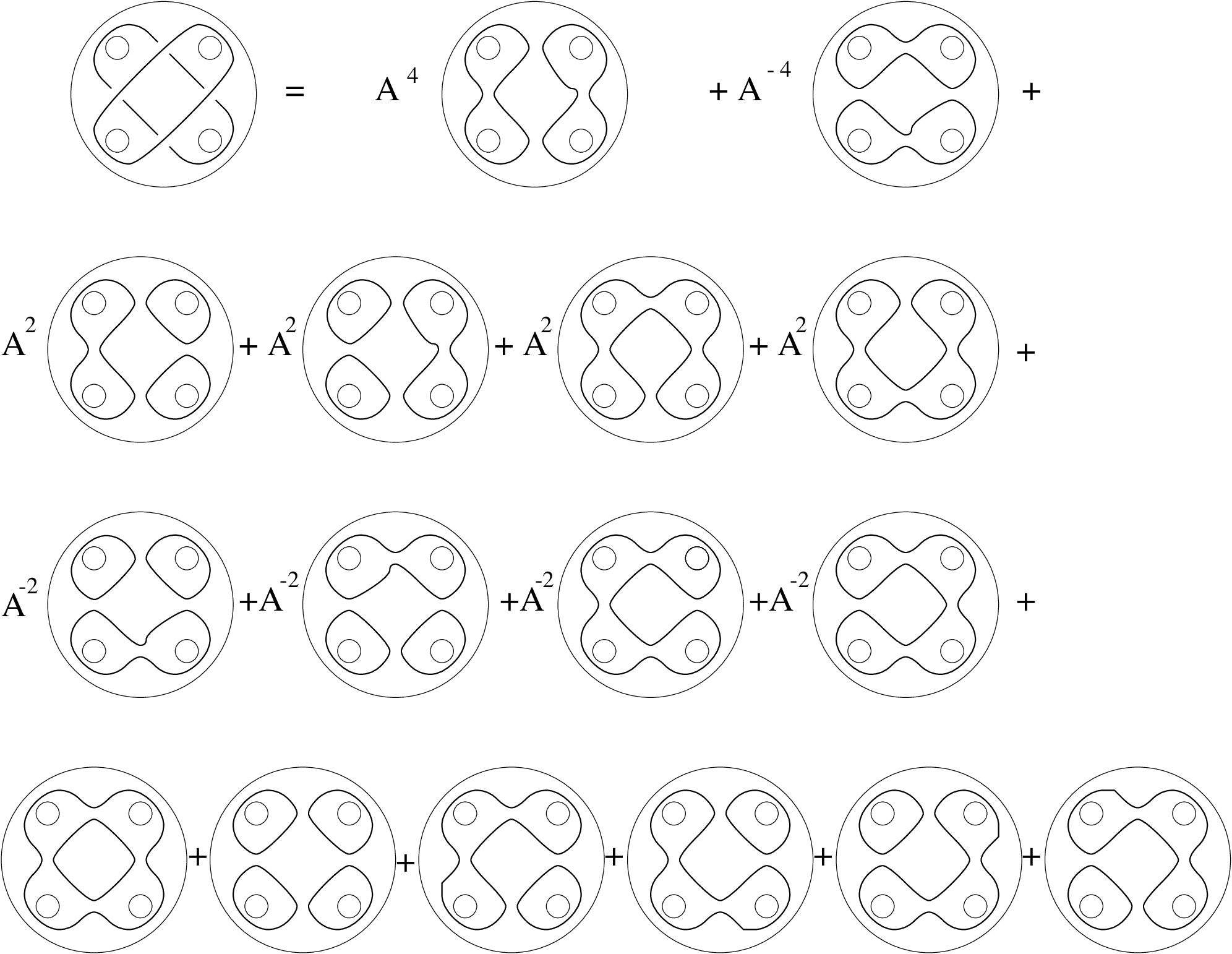}{8cm}}
 \centerline{Figure 7}
\end{figure}
 
Since all terms in the above equation, except $(A^2+A^{-2})[K_{S_0}]$ involve
only elements $[K_S],$ for $\card{S}<\card{S_0},$ we get
$(A^2+A^{-2})[K_{S_0}]\in P.$ Therefore $[K_{S_0}]\in P,$ what contradicts
our earlier assumption.

\item The minimal number of intersections between $K_{\{i_1,i_3\}}$ and 
$K_{S_0\setminus \{i_1,i_3\}}$ is three.\\
Notice that then $K_{\{i_1,i_3\}}, K_{S_0\setminus \{i_1,i_3\}}\subset
F_{g,n},$ where $g\geq 1,n\geq 3.$ We will assume, for simplicity, that 
$g=1, n=3.$
The proof for all other surfaces $F_{g,n}$ is identical.
The assumptions $g=1, n=3$ imply that $S_0=\{1,2,3,4\}.$
The product of $K_{\{i_1,i_3\}}$ and $K_{S_0\setminus \{i_1,i_3\}}$ is 
shown in Fig. 8.

\medskip
\centerline{\rysunek{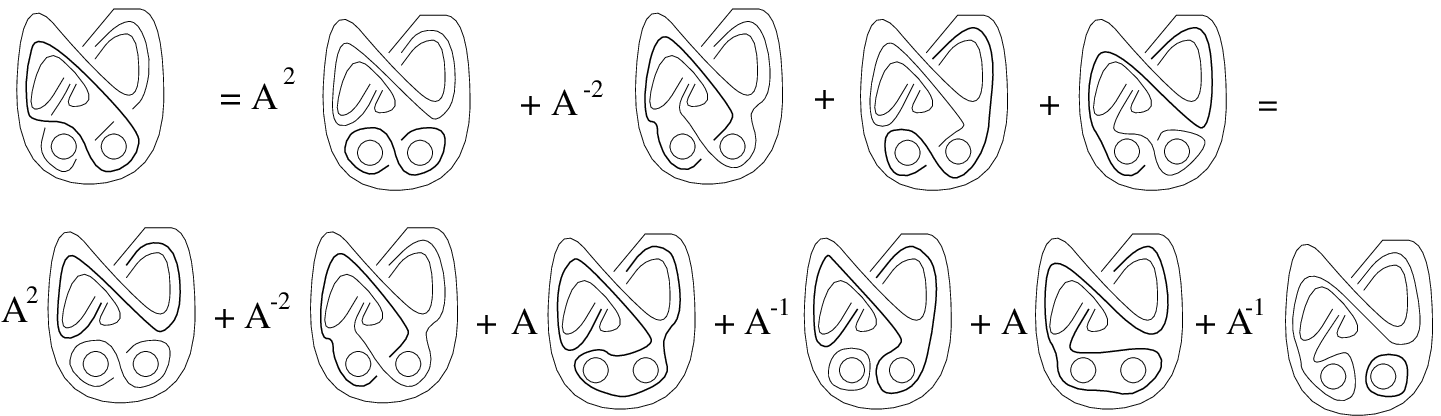}{3cm}}
\centerline{Figure 8}

Notice that any knot of the form \knot{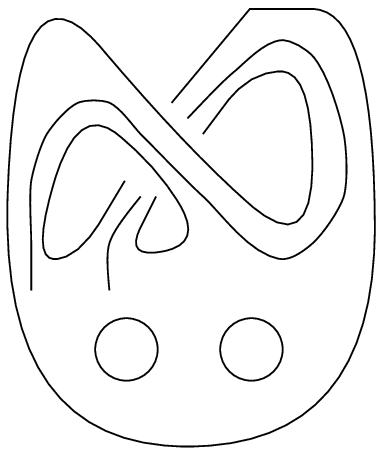} considered as an element of a
skein module is equal to $A^{-1} \knot{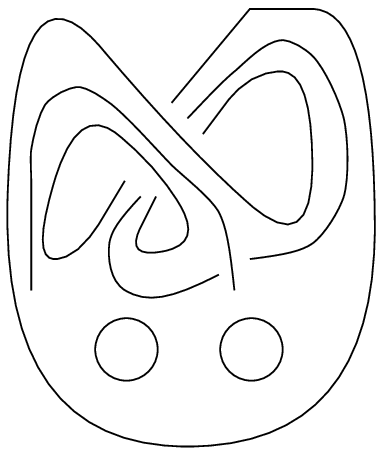}- A^{-2} \knot{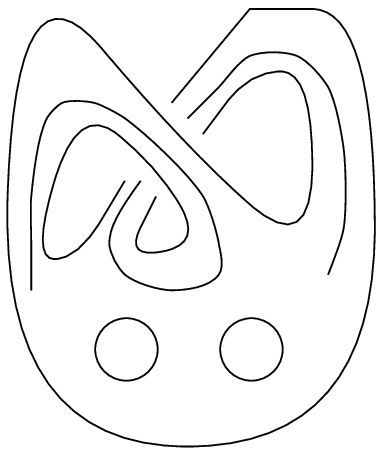}.$
Therefore $\knot{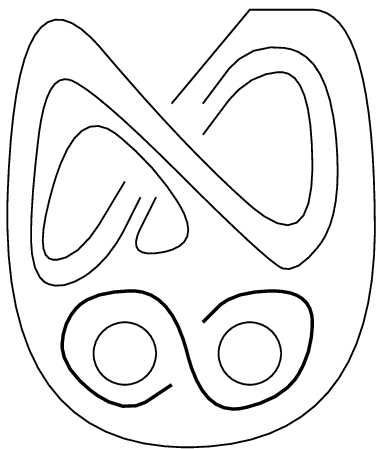}, \knot{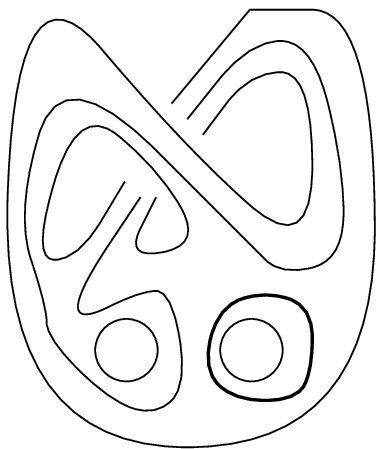} \in P.$ Similarly $\knot{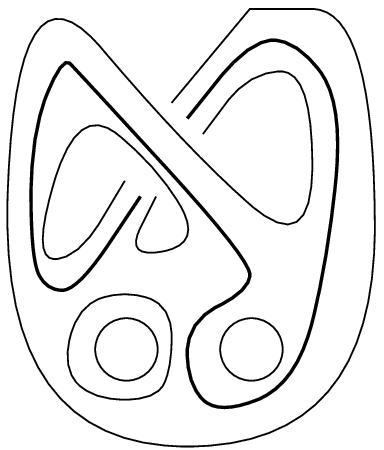}\in P.$
Moreover, $\knot{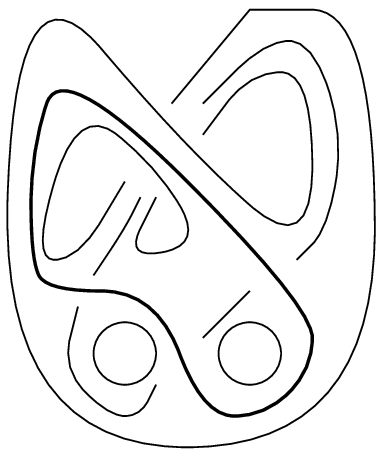}, \knot{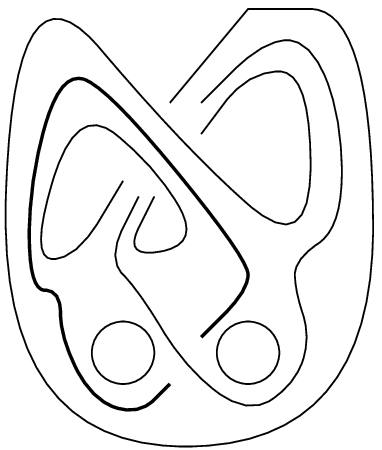}\in P.$
Hence $A \knot{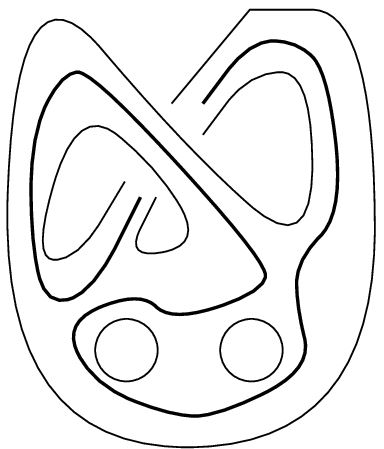} +A\knot{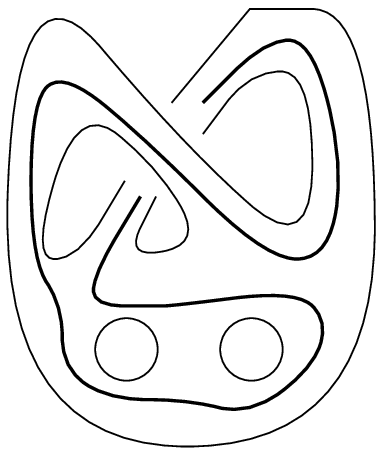}\in P.$ But $A \knot{22} +A\knot{23}= 
A^2 \knot{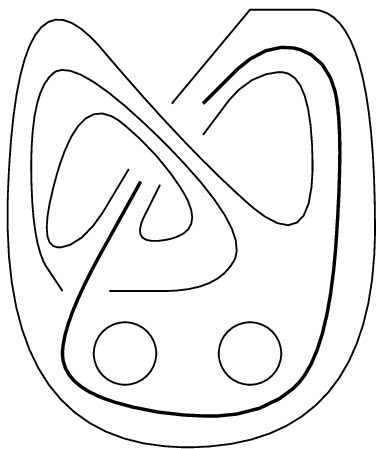}-A^3\knot{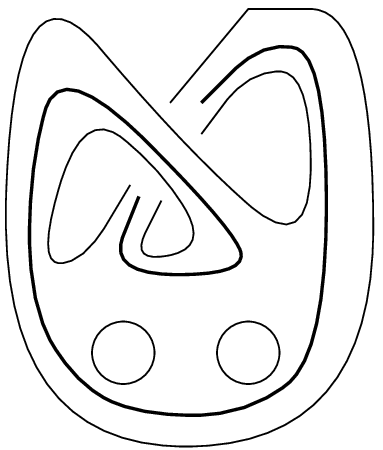}+\knot{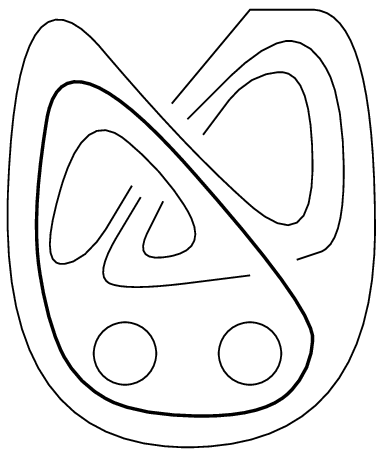}- A^{-1}\knot{24}.$
Therefore, $(A^3+A^{-1})\knot{24}\in P.$
Since $A^3+A^{-1}=A(A^2+A^{-2})$ is invertible in $R$ we get a contradiction.

\item The minimal number of intersections between $K_{\{i_1,i_3\}}$ and 
$K_{S_0\setminus \{i_1,i_3\}}$ is two.\\
Notice that then $g\geq 2.$ In this part of the proof we will assume, for 
simplicity, that $g=2$ and $n=1.$ (For all different surfaces the proof is 
identical). Then $S_0=\{1,2,3,4\}$ and $[K_{S_0}]=\knot{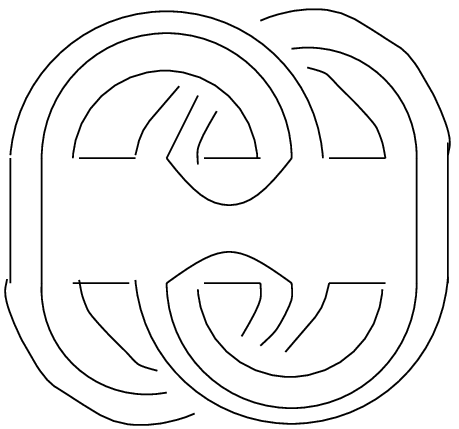}.$
Observe that $\knot{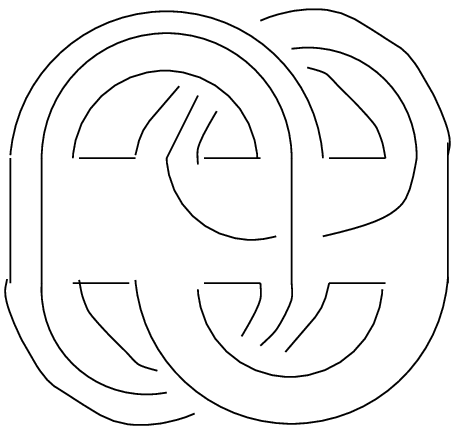} =A\cdot \knot{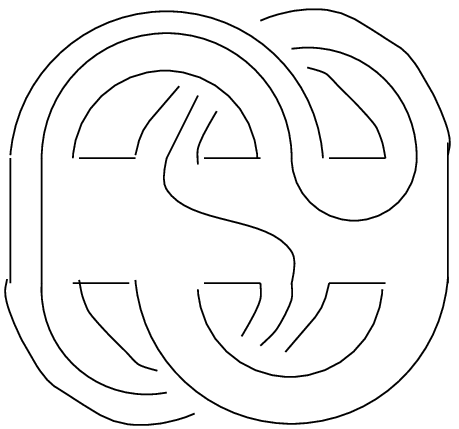} +A^{-1} \knot{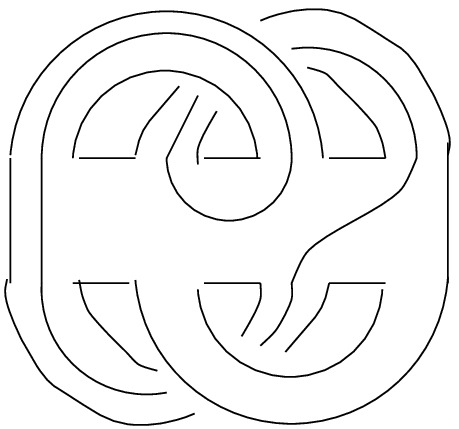}.$
Since $\knot{7},\knot{9}\in P,$ we get $\knot{8}\in P.$
Hence 
\begin{equation}
\knot{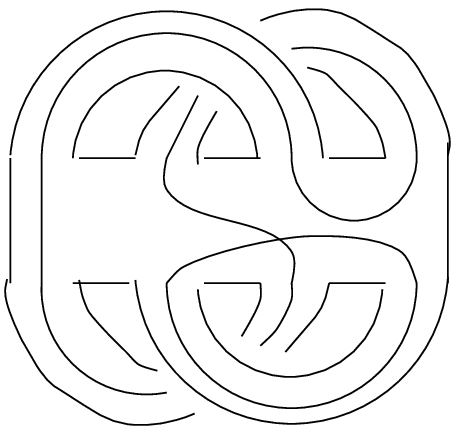}= A^{\epsilon_1}\cdot \knot{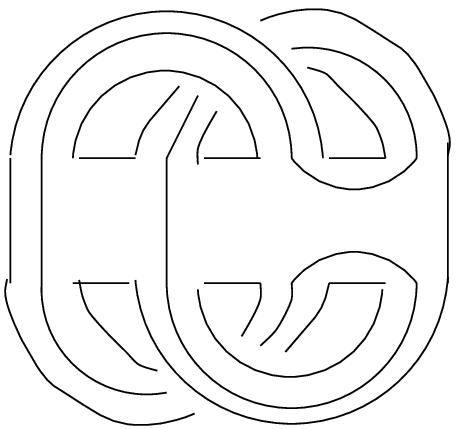}+A^{-\epsilon_1}\cdot \knot{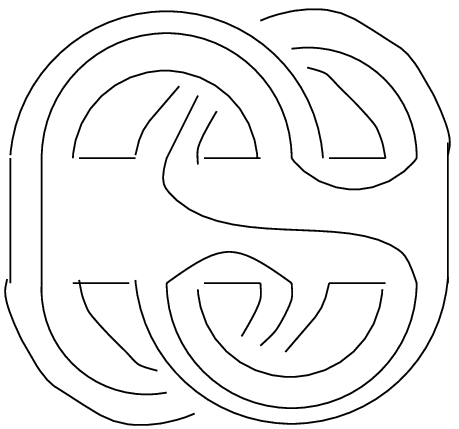}\in P,
\label{eq1}
\end{equation}
where $\epsilon_1=+1$ or $\epsilon_1=-1$ depending on the type of the crossing 
in $\knot{10}.$
Moreover,
\begin{equation}
\knot{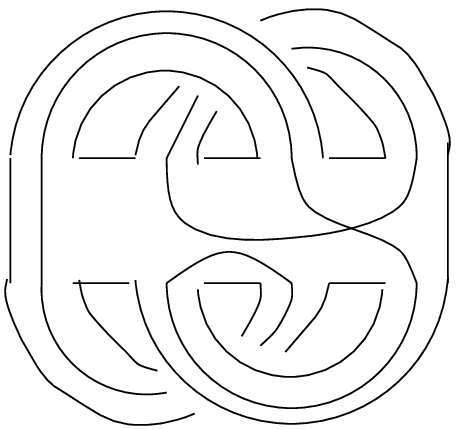}= A^{\epsilon_2}\cdot \knot{11} + A^{-\epsilon_2}\knot{12}\in P.
\label{eq2}
\end{equation}
After combining equations (\ref{eq1}) and (\ref{eq2}) we get
\begin{equation}
\knot{4}-A^{-2\epsilon_1}A^{-2\epsilon_2}\knot{12}\in P.\label{eq3}
\end{equation}
Similarly, we have
\begin{equation}
\knot{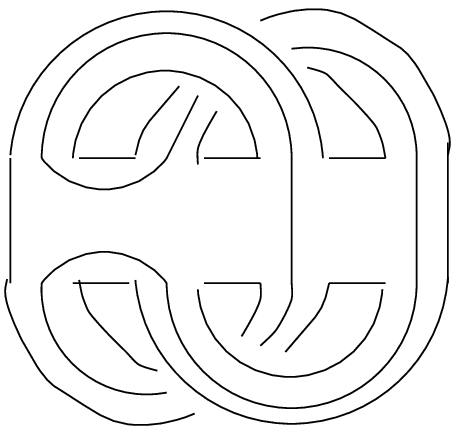}-A^{-2\epsilon_1'}A^{-2\epsilon_2'}\knot{12}\in P.\label{eq4}
\end{equation}
Finally we use the following equation:
\begin{equation}
\knot{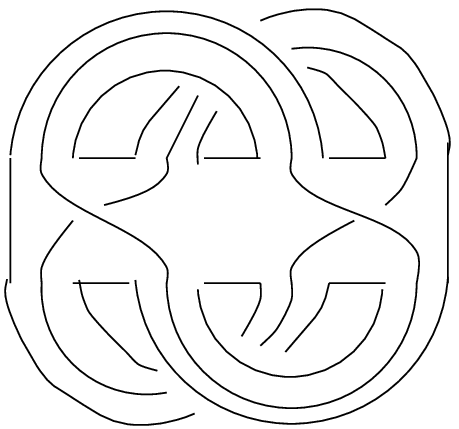}=A^2\cdot \knot{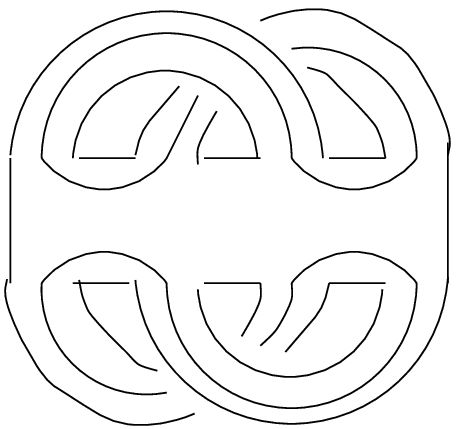} +A^{-2}\cdot \knot{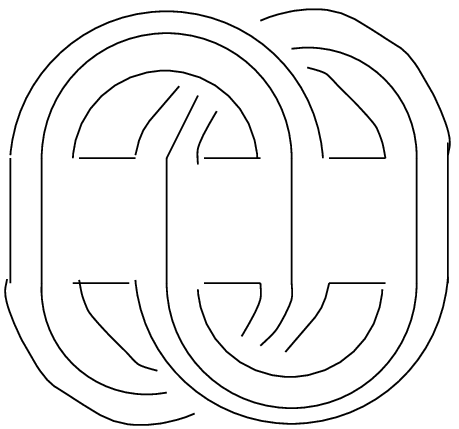}+ \knot{4}+\knot{5}.
\label{eq5}
\end{equation}
 
Since the first three links in the above equation belong to $P$ we see
that $\knot{4}+\knot{5}\in P.$
Therefore, after combining the above equation with equations (\ref{eq3}) and
(\ref{eq4}) we get $(A^{-2\epsilon_1-2\epsilon_2}+ A^{-2\epsilon_1'-
2\epsilon_2'})\cdot \knot{12}\in P,$ where $\epsilon_1,\epsilon_2,\epsilon_1',
\epsilon_2'\in \{+1,-1\}.$ In particular, we get $(A^4+1)\cdot \knot{12}\in P.$
But $A^4+1=(A^2+A^{-2})A^2$ is invertible in $R.$ Hence $[K_{S_0}]=
\knot{12}\in P.$
\end{enumerate}

This completes the proof of Theorem \ref{8.1}.\\
\end{proof}

\begin{remark}\label{8.2}\ \\
One can observe that sometimes it is possible to reduce the number of 
generators of ${\cal S}_{2,\infty}(F_{g,n}\times I;R,A)$ even if $A^2+A^{-2}$
is not invertible in $R.$ In particular, it is always possible to do this if 
at least one element of the ideal $(2,A^2+A^{-2})\triangleleft R$ is 
invertible in $R$ (i.e. $(2,A^2+A^{-2})=R$) and $g\geq 1,n\geq 2.$
\end{remark}

Theorems \ref{6.1} and \ref{8.1} may be useful for estimating Heegaard
genera of 3-manifolds. Notice that if $H_1\cup H_2=M$ is a Heegaard
splitting of $M$ then the embedding $i:H_1\to M$ induces an
epimorphism of $R$-algebras ${\cal S}_{2,\infty}(H_1;R,-1)\to 
{\cal S}_{2,\infty}(M;R,A).$ 
Therefore Theorem \ref{6.1} and Theorem \ref{8.1} imply the following fact.

\begin{proposition}\label{8.3}\ \\
Let $g(M)$ be the Heegaard genus of a compact 3-manifold $M.$ Let us denote the
minimal number of generators of any $R$-algebra $A$ by $rank(A).$ Then
\begin{enumerate}
\item $g(M)\geq log_2\biggl(\ rank\biggl({\cal S}_{2,\infty}
(M;R,-1)\biggr)+1\biggr)$
\item $g(M)+ {g(M) \choose 2}+{g(M) \choose 3}\geq rank({\cal S}_{2,\infty}
(M;R,-1))$ if ${1\over 2}\in R.$
\end{enumerate}
\end{proposition}

The above inequalities are very closely related to the obvious
inequality
\begin{equation}
g(M)\geq rank(\pi_1(M)), \label{ineq}
\end{equation}
where $rank(G)$ denotes the minimal number of generators of $G.$
However, the inequalities of Proposition \ref{8.3}
have some advantages over (\ref{ineq}). Notice that
$rank({\cal S}_{2,\infty}(M;\C,-1))\geq rank(\C[X(G)]),$ where $G=\pi_1(M).$
In practice, the minimal number of generators of $\C[X(G)]$ can be estimated 
by using methods of computational algebraic geometry. For example, one can 
calculate dimensions of tangent spaces of some irreducible components of 
$X(G)$ at singular points.
The minimal number of generators of $\C[X(G)]$ (and hence, of 
${\cal S}_{2,\infty}(M;R,-1)$) cannot be lower than any of these dimensions.

We will use this method in the latter part of this section.

\begin{theorem}\label{8.4}\ \\
Let $G$ be a group generated by $g_1,g_2,...,g_n\in G$ and let $R$ be a ring
such that $1/2\in R.$ Then
\begin{enumerate}
\item ${\cal S}(G ;R)$ can be generated by $n+{n \choose 2}+{n \choose 3}$ 
elements of the form $[g_i], [g_jg_k],$ $[g_sg_tg_v],$ for 
$1\leq i\leq n, 1\leq j<k\leq n, 1\leq s<t<v\leq n.$
\item If $G$ is abelian then
${\cal S}(G; R)$ can be generated by $n+{n \choose 2}$ elements
$[g_i], [g_jg_k], 1\leq i\leq n, 1\leq j<k\leq n.$ 
\end{enumerate}
\end{theorem}

{\it Proof:}
\begin{enumerate}
\item Let $\pi:F_n=<a_1,...,a_n>\to G$ be an epimorphism, $\pi(a_i)=g_i, 
i=1,2,...,n.$
By Fact \ref{2.6}(3) $\pi$ induces an epimorphism $\pi_*: {\cal S}(F_n; R)\to
{\cal S}(G; R)\ \pi_*([a_i])=[g_i].$ Therefore it is enough to show that the 
elements $[a_i],[a_ja_k],[a_sa_ta_v],$ for $1\leq i\leq n, 1\leq j<k
\leq n, 1\leq s<t<v\leq n,$ generate ${\cal S}(F_n; R).$ 

Since ${1\over 2}\in R,$ Theorem \ref{8.1} implies that
${\cal S}_{2,\infty}(F_{0,n+1}\times I;R,-1)$ is generated by elements 
$[K_{\{i\}}],[K_{\{j,k\}}], [K_{\{s,t,v\}}].$ By Theorem \ref{2.8}
$$\hat{\xi}:{\cal S}_{2,\infty}(F_{0,n+1}\times I;R,-1) \to {\cal S}
(<a_1,a_2,...,a_n>; R)$$ is an isomorphism of algebras such that
$\hat{\xi}([K_{\{i\}}])=-[a_i].$ Therefore  ${\cal S}(F_n; R)$ is generated
by elements $[a_i],[a_ja_k],[a_sa_ta_v],$ for $1\leq i\leq n, 
1\leq j<k\leq n, 1\leq s<t<v\leq n,$

\item Let $\Z^n$ be a free abelian group generated by $g_1',g_2',...,g_n'.$
Let $\pi:\Z^n\to G$ be an epimorphism, $\pi(g_i')=g_i.$
As in the proof of (1) it is enough to show that ${\cal S}(\Z^n; R)$ can be
generated by $g_i',g_j'g_k',$ for $1\leq i\leq n,1\leq j<k\leq n.$

Let us consider the group ring $R\Z^n\simeq 
R[x_1^{\pm 1},x_2^{\pm 1},..., x_n^{\pm 1}].$ Let\\
$R[x_1^{\pm 1},x_2^{\pm 1},..., x_n^{\pm 1}]^{sym}$ denote a subring of
$R[x_1^{\pm 1},x_2^{\pm 1},..., x_n^{\pm 1}]$ consisted of constant elements
under an involution $\tau:R[x_1^{\pm 1},x_2^{\pm 1},..., x_n^{\pm 1}]\to
R[x_1^{\pm 1},x_2^{\pm 1},..., x_n^{\pm 1}], \tau(x_i)=x_i^{-1},$
for $i=1,2,...,n.$

We have proved in [P-S-1] (Theorem 2.3) that there exists an isomorphism
${\hat \phi}:{\cal S}(G; R)\to R[x_1^{\pm 1},x_2^{\pm 1},..., x_n^{\pm 1}]
^{sym}$ such that ${\hat \phi}([g_i])=x_i+x_i^{-1}.$ Therefore the proof of 
Theorem \ref{8.4}(2) will be completed if we show the following lemma.
\end{enumerate}

\begin{lemma}\label{8.5}\ \\
The ring $R[x_1^{\pm 1},x_2^{\pm 1},..., x_n^{\pm 1}]^{sym}$ is
generated by elements $x_i+x_i^{-1},x_jx_k+x_j^{-1}x_k^{-1},$
$1\leq i\leq n, 1\leq j<k\leq n.$
\end{lemma}

\begin{proof}
Let $a_i=x_i+x_i^{-1}, b_i=x_i-x_i^{-1},$ for $i=1,2,...,n.$ Since
$1/2\in R,$ $x_i={a_i+b_i \over 2}, y_i={a_i-b_i \over 2}$ and the ring 
$R[x_1^{\pm 1},x_2^{\pm 1},..., x_n^{\pm 1}]$ is generated by 
$a_1,...,a_n,b_1,...,b_n.$ Every element $v\in
R[x_1^{\pm 1},x_2^{\pm 1},..., x_n^{\pm 1}]^{sym}$ is a sum 
$v=\sum m_i$ such that every $m_i$ is a monomial of the form
$m_i=ca_1^{\alpha_1}...a_n^{\alpha_n}b_1^{\beta_1}...b_n^{\beta_n},$ where
$c\in R, \alpha_i,\beta_i\in\{0,1,2,...\}.$ Then $v={v+\tau(v) \over 2}=
\sum {m_i+\tau(m_i) \over 2}.$ Notice that $\tau(a_i)=a_i, \tau(b_i)=-b_i.$
Hence
$${m_i+ \tau(m_i) \over 2}= \cases{m_i & if $\sum_{j=1}^n \beta _j$ is even 
\cr 0 & otherwise \cr}$$ Therefore $v$ is a sum of monomials of the form
$ca_1^{\alpha_1}...a_n^{\alpha_n}b_1^{\beta_1}...b_n^{\beta_n},$ where
$\sum_{j=1}^n \beta_j$ is even. Hence
$R[x_1^{\pm 1},x_2^{\pm 1},..., x_n^{\pm 1}]^{sym}$ is generated by elements
$a_i, b_jb_k, 1\leq i\leq n, 1\leq j<k\leq n.$
But $a_ja_k+b_jb_k=2(x_jx_k+x_j^{-1}x_k^{-1}).$
Therefore $R[x_1^{\pm 1},x_2^{\pm 1},..., x_n^{\pm 1}]^{sym}$ is
generated by $a_i=x_i+x_i^{-1}$ and $x_jx_k+x_j^{-1}x_k^{-1},$
for $1\leq i\leq n, 1\leq j<k\leq n.$
\end{proof}

Our next goal is to prove a theorem analogous to Theorem \ref{6.3}(3) but
concerning rings $R$ such that $A^2+A^{-2}$ is invertible in $R.$

\begin{theorem}\label{8.6}\ \\
Let a ring $R$ be an integral domain of characteristic $0.$
Let $G$ be a free abelian or nonabelian group on $n$ generators.
Then the numbers of generators given in Theorem \ref{8.4} are minimal i.e.
\begin{enumerate}
\item ${\cal S}(F_n; R)$ cannot have fewer then 
$n+{n \choose 2}+{n \choose 3}$ generators.
\item ${\cal S}(\Z^n; R)$ cannot have fewer then 
$n+{n \choose 2}$ generators.
\end{enumerate}
\end{theorem}

We are going to prove Theorem \ref{8.6} using methods of algebraic
geometry.
Let $G$ be any group and let $\chi_0:G\to \C$ be the trace of the trivial
$Sl_2(\C)$-representation (i.e. $\chi_0=2$). For a proof of Theorem \ref{8.6}
we will need to calculate dimensions of tangent spaces $T_{\chi_0}X(G),$
for $G=F_n,\Z^n.$ These dimensions will estimate from below the minimal
numbers of generators of the corresponding skein algebras.

Lustig and Metzler showed in [L-M] (Theorem 1) the following result.

\begin{fact}\label{8.7}\ \\
The tangent space $T_{\chi_0}X(F_n)$ has dimension $n+{n \choose 2}+
{n \choose 3}.$
\end{fact}

Using Lustig's and Metzler's idea and the lemma below we will show that\ \\
$ dim T_{\chi_0} X(\Z^n) = n + {n \choose 2}.$

\begin{lemma}\label{8.8}\ \\
Let $e_1=(1,0), e_2=(0,1), e_1,e_2\in \Z^2.$ Then any $Sl_2(\C)$-character 
$\chi$ of $\Z^2$ is uniquely determined by the values $a_1=\chi(e_1),
a_2=\chi(e_2)$ and $b=\chi(1,1).$ Therefore $X(\Z^2)\subset \C^3.$ Moreover,
$X(\Z^2)=\{(a_1,a_2,b)\in \C^3: a_1^2+a_2^2+b^2-a_1a_2b-4\}.$
\footnote{Although this is a well-known fact, we had a trouble with 
finding a good reference for a proof. For that reason we decided to give an
elementary proof here. We have seen earlier many times that many algebraic 
results concerning character varieties can be proven using topological
methods (i.e. skein algebras). The description of the skein algebra 
${\cal S}(S^1\times S^1\times I; R,A)$ given in [B-P] immediately gives a 
topological proof of Lemma \ref{8.8}.}
\end{lemma}

\begin{proof}
By Theorem \ref{8.4}(2) ${\cal S}(\Z^2; \C)$ is generated by $(1,0),(0,1),
(1,1)\in \Z^2$ and therefore, by the remarks preceding Theorem
\ref{7.3}, $\C[X(\Z^2)]$ is generated by $\tau_{(1,0)}, \tau_{(0,1)},
\tau_{(1,1)}.$ Hence we can assume that 
$X(\Z^2)\subset \C^3$ and any character $\chi\in X(\Z^2)$ has coordinates 
$(a_1,a_2,b)\in \C^3,$ where $a_1=\chi(1,0), a_2=\chi(0,1), b=\chi(1,1).$ 
Let $P(x,y,z)=x^2+y^2+z^2-xyz-4.$ We are going to show that for any character 
$\chi\in X(\Z^2)$ $P(a_1,a_2,b)=0.$

Let $\chi$ be the trace of a representation $\rho:\Z^2\to Sl_2(\C).$ Since
we can replace $\rho$ by any of its conjugate representations we can assume 
that $\rho(1,0)$ has an upper triangular form, 
$\rho(1,0)=\left(\matrix{c_1 & d_1\cr 0 & c_1^{-1} \cr}\right).$
Since $\rho(0,1)$ commutes with $\rho(1,0),$ the matrix $\rho(0,1)$ has
also an upper triangular form, 
$\rho(0,1)=\left(\matrix{c_2 & d_2\cr 0 & c_2^{-1} \cr}\right).$
Then $\chi=tr\circ \rho',$ where $\rho'$ is a diagonal representation
$$\rho'(1,0)=\left(\matrix{c_1 & 0\cr 0 & c_1^{-1} \cr}\right),\quad 
\rho'(0,1)=\left(\matrix{c_2 & 0\cr 0 & c_2^{-1} \cr}\right).$$
Hence $a_1=c_1+c_1^{-1}, a_2=c_2+c_2^{-1}, b=c_1c_2+c_1^{-1}c_2^{-1}$
and straightforward calculations show that $P(a_1,a_2,b)=0.$

Suppose now that $P(a_1,a_2,b)=0$ for some $(a_1,a_2,b)\in \C^3.$
We are going to show that $a_1=\chi(1,0), a_2=\chi(0,1), b=\chi(1,1)$
for some character $\chi\in X(\Z^2).$ Let $c_1,c_2\in \C$ be such that
$c_1+c_1^{-1}=a_1, c_2+c_2^{-1}=a_2$ and let $b'=c_1c_2+c_1^{-1}c_2^{-1},
b''=c_1c_2^{-1}+c_1^{-1}c_2.$ There are three possible cases:
\begin{enumerate}
\item $c_1,c_2\ne \pm 1.$\\
Since the points $(a_1,a_2,b'), (a_1,a_2,b'')\in \C^3$ correspond to the traces
of representations $\rho,\rho':\Z^2\to Sl_2(\C),$
$$\rho(1,0)=\left(\matrix{c_1 & 0\cr 0 & c_1^{-1} \cr}\right),\quad
\rho(0,1)=\left(\matrix{c_2 & 0\cr 0 & c_2^{-1} \cr}\right)$$
$$\rho'(1,0)=\left(\matrix{c_1 & 0\cr 0 & c_1^{-1} \cr}\right),\quad 
\rho'(0,1)=\left(\matrix{c_2^{-1} & 0\cr 0 & c_2 \cr}\right),$$
the first part of the proof implies that $P(a_1,a_2,b')=P(a_1,a_2,b'')=0.$
Moreover, $b'-b''=(c_1-c_1^{-1})(c_2-c_2^{-1})\ne 0.$
Since the equation $P(a_1,a_2,z)=0$ has at most $2$ solutions,
the equation $P(a_1,a_2,b)=0$ implies that $b=b'$ or $b=b''.$ Hence the trace 
of $\rho$ or the trace of $\rho'$ is a character of $\Z^2$ with the 
coordinates $(a_1,a_2,b).$

\item $c_1=1$ or $c_2=1.$\\
Suppose that $c_2=1$ (for $c_1=1$ the proof is identical).
Then $a_2=2$ and $a_1^2+b^2-2a_1b=0.$ Therefore $a_1=b$ and the trace of the 
representation $\rho:\Z^2\to Sl_2(\C),$
$$\rho(1,0)=\left(\matrix{c_1 & 0\cr 0 & c_1^{-1} \cr}\right), \quad
\rho(0,1)=\left(\matrix{1 & 0\cr 0 & 1 \cr}\right)$$ is the desired
$Sl_2(\C)$-character of $\Z^2.$

\item $c_1=-1$ or $c_2=-1.$\\
Suppose that $c_2=-1.$ Then $a_2=-2, a_1=-b,$ and the trace of the 
representation $\rho:\Z^2\to Sl_2(\C),$
$$\rho(1,0)=\left(\matrix{c_1 & 0\cr 0 & c_1^{-1} \cr}\right),\quad 
\rho(0,1)=\left(\matrix{-1 & 0\cr 0 & -1 \cr}\right)$$ is the desired 
$Sl_2(\C)$-character of ${\Z}^2.$

\end{enumerate}
\end{proof}

\begin{theorem}\label{8.9}\ \\
$dim T_{\chi_0}X(\Z^n) = n+{n \choose 2}={n+1 \choose 2}.$
\end{theorem}

\begin{proof}
Let $e_1=(1,0,...0),e_2=(0,1,...,0),...,e_n=(0,0,...,1)\in \Z^n.$ By Theorem 
\ref{8.4}(2) and the remarks preceding Theorem \ref{7.3}, $\C[X(\Z^n)]$
is generated by $\tau_{e_i}, \tau_{e_j+e_k},$ for $1\leq i\leq n, 1\leq 
j<k\leq n.$ Therefore we assume that $X(\Z^n)\subset \C^{n+{n \choose 2}}$ and
any $\chi\in X(\Z^n)$ has coordinates $c_i=\chi(e_i)$ and 
$d_{k,l}=\chi(e_k+e_l), k<l.$ 

Let $h_{s,t}:\Z^n\to \Z^2$ be a homomorphism $h_{s,t}(x_1,...,x_n)=(x_s,x_t),$
for $s,t\in \{1,...,n\}, s<t.$ Let $f_{s,t}:X(\Z^2)\to X(\Z^n)$ be a morphism
such that $f_{s,t}(\chi)=\chi\circ h_{s,t}.$
One can easily see that if $\chi\in X(\Z^2)$ has the coordinates
$(a_1,a_2,b)\in \C^3$ (we use the notation from the previous lemma), then
$f_{s,t}(\chi)$ has coordinates $c_i, d_{k,l},$ where
$c_i=\chi\circ h_{s,t}(e_i),d_{k,l}=\chi\circ h_{s,t}(e_k +e_l).$
Therefore
$$c_i=\cases{a_1 & if $i=s$\cr a_2 & if $i=t$ \cr 2 & otherwise \cr}
d_{k,l}=\cases{a_1 & if $s\in \{k,l\}, t\not \in \{k,l\}$\cr 
a_2 & if $t\in \{k,l\}, s\not \in \{k,l\}$\cr b & if $
\{s,t\}=\{k,l\}$\cr 2 & otherwise \cr}$$

Since all partial derivatives of $a_1^2+a_2^2+b^2-a_1a_2b-4$ disappear
for $(a_1,a_2,b)=(2,2,2),$ the variety $X(\Z^2)$ has the tangent space
at $\chi_0$ isomorphic to $\C^3.$ Therefore the linear map 
$f_{s,t*}:T_{\chi_0} X(\Z^2)\to T_{\chi_0}X(\Z^n)\subset \C^{n+{n \choose 2}}$
has a matrix presented in Figure 9. A star in an entry of the matrix denotes 
any number.

\begin{figure}
\centerline{\rysunek{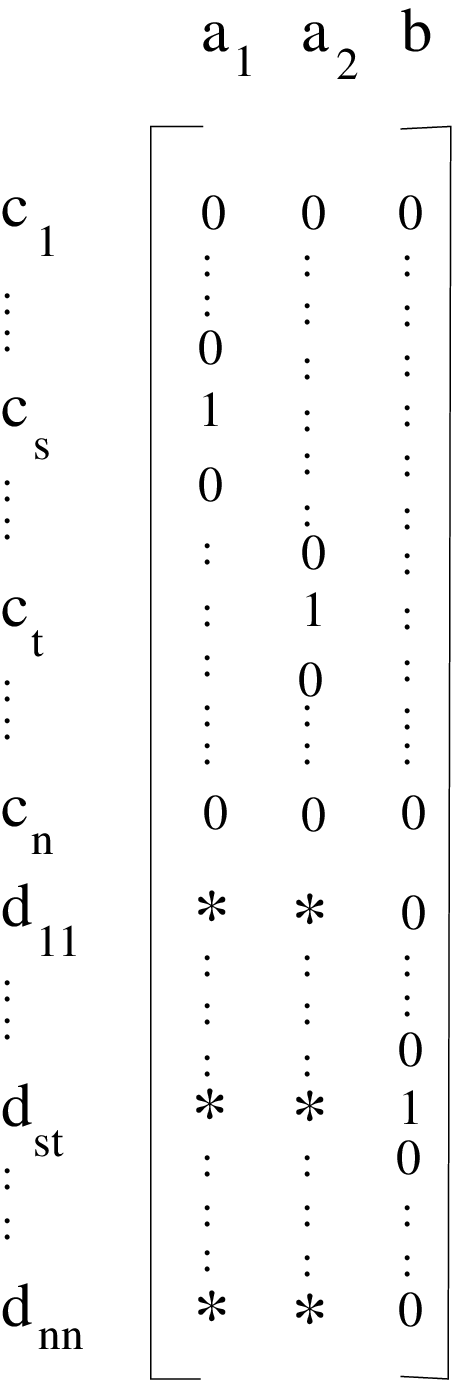}{9cm}}
\centerline{Figure 9: The matrix of $f_{s,t*}.$}
\end{figure}

One can easily see that the subspaces $Im f_{s,t*}\subset \C^{n+{n \choose 2}}$
for all $s,t\in\{1,...,n\}, s<t$ span the space $\C^{n+{n \choose 2}}.$ Since
$X(\Z^n)\subset \C^{n+{n \choose 2}},$ 
$dim T_{\chi_0} X(\Z^n)=n+{n \choose 2}.$
\end{proof}

{\it Proof of Theorem \ref{8.6}:}\ \\
Let $G=F_n$ and $N=n+{n \choose 2}+{n \choose 3}$ (respectively:
$G=\Z^n$ and $N=n+{n \choose 2}$). We denote the algebraic closure of the 
field of fractions of $R$ by $K.$ The assumptions about $R$ in Theorem
\ref{8.6} imply that $K$ has characteristic $0.$ We have noticed in the remark
preceding Theorem \ref{7.3} that there exists an epimorphism 
$\phi: {\cal S}(G; K)\to K[X(G)],$ where X(G) denotes the $Sl_2(K)$-character
variety of $G.$
Since the number of generators of $K[X(G)]$ cannot be smaller than the 
dimension of $T_{\chi_0}X(G),$ $K[X(G)],$ and therefore ${\cal S}(G; K),$ 
have at least $N$ generators
\footnote{Here we use the fact that Lustig's and Metzler's theorem
as well as our result on dimensions of tangent spaces of character varieties
and all results of Section 7
are valid for any algebraically closed field of characteristic $0.$}.
By Fact \ref{2.6}(4) ${\cal S}(G; K)={\cal S}(G; R)\otimes_R K.$ 
Hence ${\cal S}(G; R)$ has at least $N$ generators.

\section{Acknowledgments}

The authors would like to express their gratitude to D. Bullock, G. Brumfiel,
B. Goldman, C. Frohman and J. Millson for fruitful discussions. 
The second
author would like to thank V. Jones for his support during the time in
which this paper was written.

Authors addresses:

\centerline{\it Department of Mathematics, The George Washington University}
\centerline{\it 2201 G Str. room 428 Funger Hall}
\centerline{\it Washington, D.C. 20052}
\centerline{\it email: przytyck@math.gwu.edu}
\bigskip
 
\centerline{\it Department of Mathematics, University of Maryland} 
\centerline{\it College Park, MD 20742}
\centerline{\it e-mail: asikora@math.umd.edu}


\begin{thebibliography}{99}

\bibitem [Bar] {Bar} J. Barrett, Skein Spaces and Spin Structures,
University of Nottingham (1995), preprint.

\bibitem [Bau] {Bau} G. Baumslag, Topics in combinatorial group theory.
Lectures in Mathematics ETH Zurich, Birkhauser Verlag, 1993.

\bibitem [B-H]{B-H} G. W. Brumfiel, H. M. Hilden, $Sl(2)$ Representations of 
Finitely Presented Groups, {\em Contemp. Math.} {\bf 187} (1995).

\bibitem [B-1] {B-1} D. Bullock, Estimating a Skein Module with $SL_2(\C)$ 
characters, {\em Proc. Amer. Math. Soc.} {\bf 125} (1997), no. 6, 1835--1839.

\bibitem [B-2] {B-2} D. Bullock, Estimating the States of the Kauffman Bracket 
Skein Module, to appear in  Proceedings of The International Conference 
in Knot Theory, Warsaw 1995, {\em Banach Center Publications} {\bf 42}

\bibitem [B-3] {B-3} D. Bullock, A finite set of generators for the Kauffman 
bracket skein algebra, preprint (1995). 

\bibitem [B-4] {B-4} D. Bullock, Rings of $Sl_2(\C)$-characters and the 
Kauffman bracket skein module, preprint (1996).

\bibitem [B-5] {B-5} D. Bullock, The $(2,\infty)$-Skein Module 
of The Complement of a $(2,2p+1)$ Torus Knot, {\em J. Knot Theory 
Ramifications,} {\bf 4}(4) (1995), 619-632.

\bibitem [B-F-K] {B-F} D. Bullock, C. Frohman, 
J. Kania-Bartoszy\'nska, Understanding The Kauffman Bracket Skein 
Module, to appear in {\em J. Knot Theory Ramifications.}

\bibitem [B-P]{B-P} D. Bullock, J.H. Przytycki, Kauffman bracket skein module 
quantization of symmetric algebra and $so(3)$, to appear in
{\em Proc. Amer. Math. Soc.}

\bibitem [B-Z]{B-Z} G. Burde and H. Zieschang, Knots, de Gruyter, (1985).
 
\bibitem [C-S]{C-S} M. Culler and P.B. Shalen,
Varieties of group representations and splittings of 3-manifolds,
{\em Ann. of Math.} {\bf 117} (1983), 109-146.

\bibitem [C-R]{C-R} C. W. Curtis and I. Reiner,
Representation Theory of Finite Groups and Associative Algebras,
J. Wiley, New York, London (1962).

\bibitem [F-K] {F-K} R. Fricke, F. Klein, Vorlesungen \"uber die Theorie der 
automorphen Functionen, Vol. 1, pp. 365-370, B.G. Teubner, Leipzig 1897. 
Reprint: New York, Johnson reprint Corporation (Academic Press 1965).

\bibitem [Go] {Go} W. Goldman, Invariant functions on Lie groups and 
Hamiltonian flows of surface group representations, {\em Invent. Math.} 
{\bf 85} (1986), 263-302.

\bibitem [G-M] {G-M} F. Gonzalez-Acuna, J. M. Montesinos-Amilibia, 
On the character variety of group representations in $SL(2,\C)$ and 
$PSL(2,\C),$ {\em Math. Z.} {\bf 214} (1993), 627-652.
 
\bibitem [Ho] {Ho} R.D. Horowitz, Characters of free groups represented in 
two-dimensional special linear group, {\em Comm. Pure and Appl. Math.}
{\bf 25} (1972), 635-649.
 
\bibitem [H-P-1] {H-P-1} J. Hoste, J.H. Przytycki, A survey of skein modules 
of 3-manifolds; in  Knots 90, {\em Proceedings of the International Conference
on Knot Theory and Related Topics}, Osaka (Japan), August 15-19, 1990, Editor
A.~Kawauchi, Walter de Gruyter (1992), 363-379.
 
\bibitem [H-P-2] {H-P-2} J. Hoste, J.H. Przytycki, The $(2,\infty)$-skein 
module of lens spaces; A generalization of the Jones polynomial,
{\em J. Knot Theory Ramifications}, {\bf 2}(3) (1993), 321-333.
 
\bibitem [H-P-3] {H-P-3} J. Hoste, J.H. Przytycki, The $(2,\infty)$-skein 
module of Whitehead manifolds, {\em J. Knot Theory Ramifications}, {\bf 4}(3)
(1995), 411-427.

 
\bibitem [Hu]{H} T.W. Hungerford, Algebra,
Graduate Texts in Mathematics, Springer-Verlag (1974).

\bibitem [L-M] {L-M} M. Lustig, W. Metzler, Integral representations of 
$Aut F^n$ and presentation classes of groups, in Combinatorial Methods in 
Topology and Algebraic Geometry, Ed. J.R.Harper, R.Mandelbaum, 
Proceedings of a conference in honor of Arthur M.Stone, Rochester 1982, 
{\em Contemporary Mathematics}, {\bf 44} (1985), 51-67.



\bibitem [P-1] {P-1}  J.H. Przytycki, Skein modules of 3-manifolds, 
{\em Bull. Polish Acad. Science}, {\bf 39} (1-2) (1991), 91-100.

\bibitem [P-2] {P-2} J.H. Przytycki, Knot Theory. A combinatorial approach. 
Cambridge Univ. Press, to appear.


\bibitem [P-S-1] {P-S-1} J.H. Przytycki, A.S. Sikora, Skein algebra of a group,
to appear in Proceedings of The International Conference in Knot
Theory, Warsaw 1995, {\em Banach Center Publications} {\bf 42}

\bibitem [P-S-2] {P-S-2} J.H. Przytycki, A.S. Sikora, Skein algebras of 
surfaces, Preprint.

\bibitem [Ri] {Ri} R. Riley, Nonabelian representations of $2$-bridge knot 
groups, {\em Quart. J. Math Oxford Ser.} {\bf 35} (1984), 191-208.

\bibitem [S-1] {S-1} K. Saito, Representation varieties of a finitely 
generated group into $SL_2$ or $GL_2,$ RIMS Publications {\bf 958},
Kyoto University (1993).

\bibitem [S-2] {S-2} K. Saito, Character variety of representations of 
a finitely generated group in $SL_2,$ Proceedings of the 37th Taniguchi 
Symposium on Topology and Teichm\"uller Spaces, Finland 1995. 
Edited by S. Kojima {\it et al.}

\bibitem [Si] {Si} A. S. Sikora, On quaternion representations of groups, 
in preparation.

\bibitem [S-F] {S-F} A. S. Sikora, C. Frohman, A geometric method in
the theory of $SL_n$-representations of groups, preprint.

\bibitem [Th] {Th} Le Ty Kuok Tkhang\footnote{The proper spelling should be
Le Tu Quoc Thang.}, Varieties of Representations And Their 
Cohomology-Jump Subvarieties for knot groups, {\em Russian Acad. Sci. Sb. 
Math.} {\bf 78} (1994), No. 1, 187-209.

\bibitem [T-1] {T-1} V. G. Turaev, The Conway and Kauffman modules of the 
solid torus, {\em Zap. Nauchn. Sem. Lomi} {\bf 167} (1988), 79-89. 
English translation: J. Soviet Math.


\bibitem [Vo] {Vo} H. Vogt, Sur les invariants, fondamentaux des \'equations 
diff\'erentielles lin\'eaires du second ordre, {\em Ann. Sci. Ecole Norm. Sup.}
(3) 6, Suppl. 3-72 (1889) (Th\`ese, Paris).
\end{thebibliography}
\end{document}